\documentclass[a4paper,10pt]{amsart}
\usepackage{amsmath}
\usepackage{amsfonts}
\usepackage{amssymb}
\usepackage{amsthm}
\usepackage[all]{xy}
\usepackage[CJKbookmarks]{hyperref}
\usepackage{fancyhdr}
\usepackage{tikz-cd}
\usepackage{longtable}
\usepackage{array}
\usepackage{diagbox}

\makeatletter
\newcommand{\Rmnum}[1]{\expandafter\@slowromancap\romannumeral #1@}
\makeatother
\newtheorem{theorem}{Theorem}[section]
\newtheorem{lemma}[theorem]{Lemma}
\newtheorem{proposition}[theorem]{Proposition}
\newtheorem{corollary}[theorem]{Corollary}

\theoremstyle{definition}
\newtheorem{definition}[theorem]{Definition}

\newtheorem{hypothesis}[theorem]{Hypothesis}

\newtheorem{remark}[theorem]{Remark}

\numberwithin{equation}{section}
\AddToHook{cmd/theequation/before}{E}

\DeclareMathOperator{\cd}{cd}

\DeclareMathOperator{\depth}{depth}

\DeclareMathOperator{\gEnd}{\underline{End}}

\DeclareMathOperator{\gExt}{\underline{Ext}}

\DeclareMathOperator{\gr}{gr}

\DeclareMathOperator{\gldim}{gldim}

\DeclareMathOperator{\Gr}{Gr}

\DeclareMathOperator{\gHom}{\underline{Hom}}

\DeclareMathOperator{\injdim}{injdim}

\DeclareMathOperator{\im}{Im}

\DeclareMathOperator{\pdim}{pdim}

\DeclareMathOperator{\RHom}{RHom}

\DeclareMathOperator{\Tor}{Tor}

\DeclareMathOperator{\CMreg}{CMreg}
\DeclareMathOperator{\cmreg}{cmreg}
\DeclareMathOperator{\Exreg}{Ex-reg}
\DeclareMathOperator{\exreg}{ex-reg}
\DeclareMathOperator{\extreg}{extreg}
\DeclareMathOperator{\torreg}{torreg}
\DeclareMathOperator{\Torreg}{Torreg}
\DeclareMathOperator{\Extreg}{Extreg}
\DeclareMathOperator{\reg}{reg}
\DeclareMathOperator{\sdeg}{sup.deg}
\DeclareMathOperator{\ideg}{inf.deg}
\DeclareMathOperator{\ASreg}{ASreg}
\DeclareMathOperator{\asreg}{asreg}

\DeclareMathOperator{\D}{D}

\newcommand\uxm[2]{u^{#1}{({#2})}}
\newcommand\lxm[2]{l^{#1}{({#2})}}

\usepackage{makecell}

\begin{document}
\title[Weighted Numerical Homological Regularities]{Weighted Numerical Homological Regularities over positively graded algebras}

\author{Quanshui Wu}
\address{School of Mathematical Sciences \\
Fudan University \\
Shanghai, 200433 \\
 China}
\email{qswu@fudan.edu.cn}

\author{Bojuan Yi}
\address{School of Mathematical Sciences, East China Normal University, Shanghai 200241, China}
\email{bjyi@math.ecnu.edu.cn}

\thanks{This research has been supported by the NSFC (Grant No. 12471032) and the National Key Research
and Development Program of China (Grant No. 2020YFA0713200).}



\dedicatory{}

\begin{abstract}
Let $A$ be a noetherian positively graded algebra. We investigate several weighted numerical homological regularities -- weighted CM-regularity, weighted Tor-regularity (Ext-regularity), and weighted Ex-regularity -- defined in terms of the weighted supremum or infimum degrees of canonical complexes associated with graded $A$-modules. Introduced by Kirkman, Won, and Zhang, weighted regularities provide finite invariants even when their unweighted counterparts are infinite. Through parametrized degree reweightings, 
they unify classical homological invariants, including unweighted numerical regularities, $\sup$, $\inf$, depth, projective dimension, and local cohomology dimension. This framework clarifies the behavior of homological regularities under degree shifts and reveals intrinsic relationships among these invariants, particularly those involving the $a$-invariant.

We establish several inequalities and equalities relating these weighted regularities to the weighted supremum or infimum degrees of derived Hom and derived tensor complexes. These results recover their unweighted counterparts and yield new relations involving the $a$-invariant for a fixed weight. Some of these results generalize the inequalities of Kirkman--Won--Zhang between weighted CM-regularity and Tor-regularity by allowing the weight to range over the entire plane $\mathbb{R}\times\mathbb{R}$. We also give conditions under which these inequalities become equalities, thereby answering a question posed by Kirkman--Won--Zhang. We use weighted AS-regularity to study the Artin–Schelter regularity of $A$. Under mild conditions, we prove that the $\xi$-AS-regularity of $A$ vanishes if and only if $A$ is an $\mathbb{N}$-graded AS-regular algebra. This generalizes an unweighted result of the present authors and a weighted result of Kirkman, Won, and Zhang. 
\end{abstract}

\maketitle

\section{Introduction}
Castelnuovo–Mumford regularity (CM-regularity for short) originated in Mumford's vanishing theorem \cite{M} for coherent sheaves $\mathcal{M}$ on the projective scheme $\mathbb{P}^n$ and was subsequently developed by Ooishi \cite{Ooi} and Eisenbud--Goto \cite{EG} into a homological invariant for graded modules over commutative graded algebras. In this setting, it is closely related to Tor-regularity, linear free resolutions, and the degrees of generators of syzygies.
In the noncommutative setting, J\o rgensen \cite{Jo2, Jo4} extended CM-regularity to connected graded algebras and established fundamental inequalities relating CM-regularity to Tor-regularity. Subsequently, R\"{o}mer \cite{R\"{o}}, Dong–Wu \cite{DW}, and Kirkman–Won–Zhang \cite{KWZ1} further refined the theory. In particular, Kirkman–Won–Zhang \cite{KWZ1} introduced numerical AS-regularity to characterize Artin–Schelter regular algebras in the noetherian connected graded case.

Building on earlier studies of the connected graded case, we investigated
numerical homological regularities of complexes over noncommutative noetherian locally finite $\mathbb{N}$-graded algebras in \cite{WY}, including CM-regularity, Tor-regularity, Ext-regularity, Ex-regularity, and their companions---the lowercase-named regularities. In \cite{WY}, we established various inequalities and equalities relating these regularities in the context of derived Hom and tensor functors. We also introduced the notions of CM-regularity homogeneous and ex-regularity homogeneous, which provide criteria under which the inequalities become equalities. Under mild hypotheses, the vanishing of numerical AS-regularity (that is, $\operatorname{ASreg}(A) = 0$) characterizes $\mathbb{N}$-graded Artin–Schelter regular algebras, extending earlier results from the connected graded setting to the case in which $A_0$ is semisimple. Likewise, $\asreg(A)=0$ characterizes finite-dimensional $\mathbb{N}$-graded algebras.

Weighted versions of classical numerical regularities, including CM-regularity, Tor-regularity, Ext-regularity, and AS-regularity, were introduced for noetherian connected graded algebras $A$ by Kirkman-Won-Zhang in \cite{KWZ2}. These weighted regularities generalize the corresponding unweighted invariants. Other invariants of a complex, including \(\sup\), \(\inf\), projective dimension, depth, and the $a$-invariant, can also be studied within the framework of weighted regularities of complexes. 
The $a$-invariant is a classical numerical invariant of graded Gorenstein algebras. In the commutative setting, it was introduced by Goto and Watanabe and is treated systematically in the monograph by Bruns and Herzog \cite{GW,BH}. The family of $a_i$-invariants generalizes the $a$-invariant by measuring the highest non‑vanishing degree of each graded local cohomology module. These invariants were studied by Trung and his collaborators; see \cite{Tr1,Tr2,HT,CHT,TW} for original results and \cite{Tr3} for a survey. For noncommutative noetherian connected graded AS‑Gorenstein algebras, the analogue of the $a$-invariant, often called the AS‑index, first appeared in the work of Artin and Schelter and encodes the degree shift of the balanced dualizing complex. Comprehensive surveys can be found in \cite{AS,Jo3,Z,JZ}. As observed in \cite{KWZ2}, the $a$-invariant can be realized as a special case of the weighted numerical‑regularity formalism.
In addition, Kirkman-Won-Zhang established fundamental inequalities relating these invariants \cite{KWZ2}. For any noetherian connected graded algebra $A$ with a balanced dualizing complex and any $X\in \mathbf{D}^b(\gr{A})$, the inequalities $\Torreg_\xi(X) \leqslant \CMreg_\xi(X) + \Torreg_\xi(k)$ and $\CMreg_\xi(X) \leqslant \Torreg_\xi(X) + \CMreg_\xi(A)$ hold for all $\xi\in \mathbb{R}$. When $X$ has finite projective dimension, the second inequality is an equality for $0 \leqslant \xi \leqslant 1$ or $\xi\ll 0$; whether it remains an equality for all $\xi\leqslant 1$ is unknown. This equality unifies several classical homological identities, including the Auslander–Buchsbaum formula. Kirkman-Won-Zhang \cite{KWZ2} also proved that $\ASreg_\xi(A)=0$ characterizes Artin–Schelter regular algebras for $\xi\leqslant 1$ and used the weighted theory to establish the finiteness of weighted Tor-regularity under finite maps and the Koszul property of sufficiently high Veronese subrings.
 
In this paper, we study weighted versions of the aforementioned numerical regularities from \cite{WY} for $\mathbb{N}$-graded algebras and resolve a question raised in \cite[Remark 3.11(2)]{KWZ2}. Throughout, all $\mathbb{N}$-graded algebras are assumed to be locally finite over a field $k$. 

Classical numerical regularities can be infinite for many noncommutative graded algebras and do not provide a unified framework for classical homological invariants of complexes (supremum, infimum, projective dimension, depth, and $a$-invariant) and identities such as the Auslander–Buchsbaum formula. Weighted regularities overcome these limitations by introducing adjustable weight parameters. The weighted framework includes all unweighted results as special cases, yields finite homological invariants, establishes a general theory, and provides new equivalent characterizations of Artin–Schelter regular algebras.
The weighted numerical regularities studied in this paper, including weighted CM-regularity, weighted Tor-regularity, weighted Ext-regularity, and weighted Ex-regularity (see Definitions \ref{CM-reg-def}), are defined in terms of the weighted supremum or infimum degrees of certain canonically associated complexes.  

Let $\xi=(\xi_0, \xi_1)\in \mathbb{R}\times \mathbb{R}$. The {\it $\xi$-supremum degree} and {\it $\xi$-infimum degree} of $X \in \mathbf{D}(\Gr{A})$ are defined respectively as 
$$\sdeg_{\xi}(X)=\sup\{ \xi_1i + \xi_0 j \mid H^i(X)_j\neq 0,i,j\in \mathbb{Z} \},$$
$$ \ideg_{\xi}(X)=\inf\{ \xi_1 i+ \xi_0 j \mid H^i(X)_j\neq 0, i,j\in \mathbb{Z} \}.$$

\begin{definition} (See Definitions \ref{CM-reg-def})\label{weighted-reg}
Suppose that $A$ is an $\mathbb{N}$-graded algebra and $S=A/J$, where $J$ is the graded Jacobson radical of $A$. Let $\xi=(\xi_0,\xi_1)\in \mathbb{R}\times \mathbb{R}$.
\begin{itemize}
    \item [(1)] For any $X\in \mathbf{D}^+(\Gr A)$, the $\xi$-CM-regularity and $\xi$-Ex-regularity $X$ are defined respectively as
$$\CMreg_{\xi}(X):=\sdeg_{\xi} (R\Gamma_A(X)), \, \Exreg_{\xi}(X):=\sdeg_{\xi} (R\gHom_A(S, X)).$$

   \item [(2)] For any $X\in \mathbf{D}^-(\Gr A)$, the $\xi$-Tor-regularity and $\xi$-Ext-regularity of $X$ are defined respectively as
$$\Torreg_{\xi}(X):=\sdeg_{\xi} (S \,\, {}^L\!\otimes_A X), \, \Extreg_{\xi}(X):=-\ideg_{\xi} (R\gHom_A(X,S)).$$
\end{itemize}
\end{definition}

For the weights $\xi=(1,1)$ and $\xi=(-1,-1)$, the weighted regularities defined above recover the unweighted versions investigated in \cite{WY}. For the weight $\xi=(1,0)$, the weighted CM-regularity defined above coincides with the $a$-invariant in \cite{GW,Tr3}, as discussed in Remarks \ref{relation-nonweighted-regularity} and \ref{a and b-invariant}.

Weighted CM-regularity is a key numerical invariant. We show that its finiteness with respect to various weights ensures that $A$ admits a balanced dualizing complex and provides a characterization when $A$ is finite-dimensional; see Propositions \ref{bdc and CM} and \ref{finite-dim alg}.

In the latter part of Section 3, we discuss the relationships among the four types of weighted regularity defined in Definition \ref{weighted-reg}. The following two theorems, the main results of Section 3, establish the fundamental relationships among these weighted regularities. The first generalizes \cite[Lemma 5.2]{DW} and \cite[Proposition 3.14]{WY}, while the second generalizes \cite[Theorem 3.16]{WY} and \cite[Theorem 3.9]{Ngu} using Theorem \ref{CM and Ext main result}.

\begin{theorem}[Proposition \ref{CM and Ext}]\label{CM and Ext main result}
Let $A$ be a noetherian $\mathbb{N}$-graded algebra, and let $\xi=(\xi_0,\xi_1)\in \mathbb{R}_{\leqslant 0}\times \mathbb{R}$. For any $0\neq X \in \mathbf{D}^-(\gr A),$
$$\Torreg_{\xi}(X)=\Extreg_{\xi}(X)\geqslant \sdeg_{\xi}(X).$$
Furthermore, if $A_0$ is semisimple and $\xi_1\geqslant\xi_0$, then 
$$\Torreg_{\xi}(X)=\sdeg_{\xi}(X).$$
\end{theorem}
 
\begin{theorem}[Theorem \ref{ex}]\label{ex-main}
Suppose that $A$ is a noetherian $\mathbb{N}$-graded algebra and $\xi=(\xi_0,\xi_1)\in \mathbb{R}_{\geqslant 0}\times \mathbb{R}$. 
Then, for any $0\neq X \in \mathbf{D}^+_{fg}(\Gr A)$, 
$$\CMreg_{\xi}(X)\leqslant \Exreg_{\xi}(X).$$
 Furthermore, if $\xi_1\leqslant\xi_0$, $A_0$ is semisimple, and $A$ has a balanced dualizing complex $R$, then 
 $$\CMreg_{\xi}(X)=\Exreg_{\xi}(X).$$
\end{theorem}

We establish several inequalities and equalities concerning the weighted supremum or infimum degrees of $R\gHom_A(X,Y)$ and $Y\, {}^L\!\otimes_A X$ for $X, Y \in \mathbf{D}(\Gr A)$ over noetherian $\mathbb{N}$-graded algebras. These results, stated in Propositions \ref{sdeg-ideg inequality 1}, \ref{sdeg-ideg inequality 11}, \ref{sdeg-ideg equality 1}, \ref{sdeg-ideg equality 2}, \ref{sdeg-ideg equality 3}, \ref{sdeg-ideg equality 4}, and \ref{sdeg-ideg equality 8}, are obtained using cohomological spectral sequences and Lemmas \ref{pre for equality 1} and \ref{pre for equality 2}. For example,
\begin{proposition} \label{sdeg-ideg inequality main}
    Let $A$ be a noetherian $\mathbb{N}$-graded algebra. Then the following hold. 
    \begin{itemize}
        \item [(1)] For any $0\neq X\in \mathbf{D}^-(\gr A)$, $0\neq Y\in \mathbf{D}^+(\Gr A)$ and $\xi=(\xi_0,\xi_1)\in \mathbb{R}_{\geqslant 0}\times \mathbb{R}$, 
   \begin{align*} 
   -\ideg_{\xi}(R\gHom_A(X,Y))\leqslant \Extreg_{\xi}(X)-\ideg_{\xi}(Y).
   \end{align*}
   \item [(2)] If $A$ has a balanced dualizing complex, then for any $0\neq X \in \mathbf{D}^-_{fg}(\Gr A)$, $0\neq Y \in \mathbf{D}^-(\gr A^o)$ with torsion cohomologies and $\xi=(\xi_0,\xi_1)\in \mathbb{R}\times \mathbb{R}$,
\begin{align*}  
   -\ideg_{\xi}(R\gHom_A(X,\D(Y)))\leqslant \CMreg_{\xi}(X)+\Extreg_{\xi}(Y).
   \end{align*}
   \item [(3)] For any $0\neq X\in \mathbf{D}^-(\gr A)$, $0\neq Y\in \mathbf{D}^-(\Gr A^{o})$ and $\xi=(\xi_0,\xi_1)\in \mathbb{R}\times \mathbb{R}$,
    \begin{align*}
    -\ideg_{\xi}(Y\, {}^L\!\otimes_A X) \leqslant \Extreg_{-\xi}(X)-\ideg_{\xi}(Y).
    \end{align*}
    \end{itemize} 
    \end{proposition}
\begin{proposition}[Proposition \ref{sdeg-ideg equality 1}, \ref{sdeg-ideg equality 3}]\label{sdeg-ideg equality main}
   Suppose $A$ is a noetherian $\mathbb{N}$-graded algebra with $A_0$ semisimple and $\xi_1\leqslant \xi_0$, where $\xi=(\xi_0,\xi_1)\in \mathbb{R}_{> 0}\times \mathbb{R}$. For any $0\neq Y \in \mathbf{D}^b(\gr {A^o})$, 
   \begin{itemize}
       \item [(1)] $\ideg_{\xi}(Y\, {}^L\!\otimes_A X) =  \ideg_{\xi}(X)+\ideg_{\xi}(Y)$ for all $0\neq X\in \mathbf{D}^b(\gr A)$ if and only if $\ideg_{\xi}(Ye_i)=\ideg_{\xi}(Y)$ for all $1\leqslant i\leqslant n$.
    \item [(2)] $-\ideg_{\xi}(R\gHom_A(X,Y))=\Extreg_{\xi}(X)-\ideg_{\xi}(Y)$ for all $0\neq X\in \mathbf{D}^b(\gr A)$ with $\pdim(X)<\infty$ if and only if  
    $\ideg_{\xi}(e_iY)=\ideg_{\xi}(Y)$ for all $1\leqslant i\leqslant n$.
   \end{itemize}
\end{proposition}
Following the approach of \cite{KWZ2, WY}, we study weighted Artin-Schelter regularity for $\mathbb{N}$-graded algebras. The $\xi$-AS regularity of a noetherian $\mathbb{N}$-graded algebra $A$ is defined by $$\ASreg_{\xi}(A) =\CMreg_{\xi}(A) + \Torreg_{\xi}(S).$$

For $\xi=(1,1)$ or $(-1,-1)$, $\ASreg_{(1,1)}(A)=\ASreg(A)$ and $\ASreg_{(-1,-1)}(A)=\asreg(A)$, where $\ASreg(A)$ and $\asreg(A)$ are defined in \cite[Definition 5.1 and 5.3]{WY}.

Parts (1) and (2) of Theorem \ref{relation between CM and Tor main} below generalize \cite[Theorems 3.3 and 3.5]{KWZ2}, and the condition on $\xi$ can be relaxed to allow $\xi$ to range over the entire plane $\mathbb{R}\times \mathbb{R}$. When $\xi=(1,1)$ or $(-1,-1)$, the theorem recovers \cite[Theorem 1.8]{WY}; when $\xi=(1,0)$, it yields new relationships between the $a$-invariant and the $b$-invariant; see Corollary \ref{relation-a-b}.
\begin{theorem}[Theorem \ref{relation between CM and Tor}]\label{relation between CM and Tor main}
 Suppose $A$ is a noetherian $\mathbb{N}$-graded algebra with a balanced dualizing complex $R$. Let $\xi=(\xi_0,\xi_1)\in \mathbb{R}\times \mathbb{R}$. Then
   \begin{itemize}
       \item [(1)] $\CMreg_{\xi}(X)\leqslant \Torreg_{\xi}(X)+\CMreg_{\xi}(A)$ for all $0\neq X \in \mathbf{D}^-(\gr A)$.
       \item [(2)] $\Torreg_{\xi}(X)\leqslant \CMreg_{\xi}(X)+\Torreg_{\xi}(S)$ for all $0\neq X \in \mathbf{D}^-(\gr A)$.
       \item [(3)] $\ASreg_{\xi}(A)\geqslant 0$.
   \end{itemize}
\end{theorem}
For $\mathbb{N}$-graded algebras, we introduce the weighted CM-regularity homogeneous and weighted Ex-regularity homogeneous properties to study the relationships among the previously defined weighted regularities. We determine when the inequalities in Theorem \ref{relation between CM and Tor main} become equalities, generalize \cite[Theorem 3.10]{KWZ2} to $\mathbb{N}$-graded algebras, and answer a question posed by Kirkman-Won-Zhang in \cite[Remark 3.11 (2)]{KWZ2}.

\begin{definition}[Definition \ref{CM-homogeneous-def}]\label{CM-homogeneous-def main} Suppose $A$ is an $\mathbb{N}$-graded algebra and $\xi=(\xi_0,\xi_1)\in \mathbb{R}\times \mathbb{R}$.
 \begin{itemize}
     \item [(1)] $A$ is called {\it  left (resp. right) $\xi$-CM-regularity homogeneous} if $\CMreg_{\xi}(A)=\CMreg_{\xi}(Ae_i)$  (resp. $\CMreg_{\xi}(A)=\CMreg_{\xi}(e_iA)$) for all $1 \leqslant i \leqslant n$.
     \item [(2)] $A$ is called {\it left (resp. right) $\xi$-Ex-regularity homogeneous}  if  $\Exreg_{\xi}({}_{A}A)=\Exreg_{\xi}(Ae_i)$ (resp. 
 $\Exreg_{\xi}(A_A)=\Exreg_{\xi}(e_iA)$) for all $1 \leqslant i \leqslant n$.
    \end{itemize}
    \end{definition} 
Every connected graded algebra $A$ is $\xi$-CM-regularity homogeneous and $\xi$-Ex-regularity homogeneous for all $\xi=(\xi_0,\xi_1)\in \mathbb{R}\times \mathbb{R}$.  

Theorem \ref{relation between CM and Tor 2 main} extends \cite[Theorem 3.10]{KWZ2} and \cite[Theorem 5.24]{WY}. It follows directly from Proposition \ref{sdeg-ideg equality 3} upon taking $Y=R$, where $R$ is the balanced dualizing complex of $A$. Moreover, if $A$ is connected graded, then 
Theorem \ref{relation between CM and Tor 2 main} implies that, for any $\xi=(\xi_0,\xi_1)\in \mathbb{R}_{>0}\times \mathbb{R}$ satisfying $\xi_1\leqslant \xi_0$, the equality $\CMreg_{\xi}(X)=\Torreg_{\xi}(X)+\CMreg_{\xi}(A)$ holds for all $0\neq X\in \mathbf{D}^b(\gr A)$ with finite projective dimension. This answers a question posed by Kirkman-Won-Zhang in \cite[Remark 3.11 (2)]{KWZ2}. 

\begin{theorem}[Theorem \ref{relation between CM and Tor 2}]\label{relation between CM and Tor 2 main}
Suppose $A$ is a noetherian $\mathbb{N}$-graded algebra with $A_0$ semisimple and $A$ admits a balanced dualizing complex $R$. Let $\xi=(\xi_0,\xi_1)\in \mathbb{R}_{>0}\times \mathbb{R}$ and $\xi_1\leqslant \xi_0$. Then the following conditions are equivalent.
   \begin{itemize}
       \item [(1)] $A$ is left $\xi$-CM-regularity homogeneous.
       \item [(2)] $\CMreg_{\xi}(X)=\Torreg_{\xi}(X)+\CMreg_{\xi}(A)$ for all $0\neq X\in \mathbf{D}^b(\gr A)$ with finite projective dimension.
   \end{itemize}
\end{theorem}
For $\xi=(1,0)$, the above theorem establishes several properties of the $a$-invariant; see Corollary \ref{relation of a-invariant}.

Building on the ideas of \cite{KWZ2} and \cite{WY}, we study the relationship between weighted AS-regularities and the $\mathbb{N}$-graded AS-regular property for noetherian $\mathbb{N}$-graded algebras, thereby generalizing results for the connected graded case \cite[Theorem 0.6]{KWZ2} and for the non-weighted case of $\mathbb{N}$-graded algebras \cite[Theorem 1.16, 1.17]{WY}.
That $A$ satisfies the Auslander–Buchsbaum formula is a key hypothesis in the proof of the following theorem. The Auslander-Buchsbaum formula holds for connected graded algebras satisfying $\chi$ \cite[Theorem 3.2]{Jo1}. Lemma \ref{AB formula} examined, using homogeneous left weighted CM-regularities, the conditions under which a noetherian $\mathbb{N}$-graded algebra satisfies the Auslander-Buchsbaum formula.

In Theorem \ref{AGP}, we apply the technique of \cite[Theorem 4.7]{IKU} to adjust the Gorenstein parameters of $\mathbb{N}$-graded AS-Gorenstein algebras. Average Gorenstein parameters were also introduced in \cite{IKU}. Using Theorem \ref{AGP}, we prove

 \begin{theorem}[Theorem \ref{AS regular 11}, \ref{AS regular 22}]\label{AS regular 11 main}
Suppose that $A$ is a ring-indecomposable, basic, noetherian locally finite $\mathbb{N}$-graded algebra with a balanced dualizing complex and that $A_0$ is semisimple. Then the following statements are equivalent.
\begin{itemize}
    \item [(1)] $A$ is $\mathbb{N}$-graded AS-regular of dimension $d$ with the average Gorenstein parameter $\ell_{av}^A \in \mathbb{Z}$.
    \item [(2)] For some integers $p_1, p_2,\dots, p_n$,  $B: =\gEnd_A(\mathop{\bigoplus}\limits_{i=1}^n Ae_{i}(p_i))$ is an $\mathbb{N}$-graded CM-algebra of dimension $d$ with $B_0$ semisimple  such that $\ASreg(B)=0$.    \item [(3)] For some integers $p_1, p_2,\dots, p_n$,  $B: =\gEnd_A(\mathop{\bigoplus}\limits_{i=1}^n Ae_{i}(p_i))$ is an $\mathbb{N}$-graded CM-algebra of dimension $d$ with $B_0$ semisimple  such that $\ASreg_{\xi}(B)=0$ for some $\xi=(\xi_0,\xi_1)\in \mathbb{R}_{>0}\times \mathbb{R}$.
    \item[(4)] For some integers $p_1, p_2,\dots, p_n$,  $B: =\gEnd_A(\mathop{\bigoplus}\limits_{i=1}^n Ae_{i}(p_i))$ is an $\mathbb{N}$-graded algebra with $B_0$ semisimple such that  $\ASreg_{\xi}(B)=0$ for some $\xi=(\xi_0,\xi_1)\in \mathbb{R}_{>0}\times \mathbb{R}$ with $\xi_1\ll 0$.
    \item[(5)] For some integers $p_1, p_2,\dots, p_n$,  $B: =\gEnd_A(\mathop{\bigoplus}\limits_{i=1}^n Ae_{i}(p_i))$ is an $\mathbb{N}$-graded algebra with $B_0$ semisimple such that  $\ASreg_{\xi}(B)=0$ for all $\xi=(\xi_0,\xi_1)\in \mathbb{R}_{>0}\times \mathbb{R}$ with $\xi_1\leqslant \xi_0$.
\end{itemize}
\end{theorem}

Lemma \ref{same gorenstein parameter} characterizes when the Gorenstein parameters of any $\mathbb{N}$-graded AS-Gorenstein algebra $A$ coincide. Consequently, the following theorem extends \cite[Theorem 0.6]{KWZ2} and \cite[Theorem 6.10]{WY} to $\mathbb{N}$-graded algebras.

\begin{theorem}[Theorem \ref{AS regular 33}, \ref{weighted AS regular 11}]\label{AS regular 33 main} 
   Suppose that $A$ is a noetherian, locally finite $\mathbb{N}$-graded algebra with a balanced dualizing complex, and that $A_0$ is semisimple. Then the following statements are equivalent.
   \begin{itemize}
       \item [(1)] $A$ is $\mathbb{N}$-graded AS-regular of dimension $d$ with the same Gorenstein parameters.
        \item [(2)] $\ASreg_{\xi}(A)=0$ for all $\xi=(\xi_0,\xi_1)\in \mathbb{R}_{>0}\times \mathbb{R}$ with $\xi_1\leqslant \xi_0$.
       \item [(3)] $\ASreg_{\xi}(A)=0$ for some $\xi=(\xi_0,\xi_1)\in \mathbb{R}_{>0}\times \mathbb{R}$ with $\xi_1\ll 0$.  
   \end{itemize}
   If furthermore, $A$ satisfies both the left and right Auslander–Buchsbaum formulas, then (1) is equivalent to 
   \begin{itemize}
       \item [(4)] $\ASreg_{\xi}(A)=0$ for some $\xi=(\xi_0,\xi_1)\in \mathbb{R}_{>0}\times \mathbb{R}$ with $\xi_1\leqslant \xi_0$.
   \end{itemize}
\end{theorem}

\section{Preliminaries} 
Let $k$ be a field. An $\mathbb{N}$-graded $k$-algebra $A=\mathop{\bigoplus}\limits_{i=0}^{\infty}A_{i}$ is called {\it locally finite} if $\dim_{k}A_{i}<\infty$ for all $i\geqslant 0$. An $\mathbb{N}$-graded $k$-algebra $A$ is called {\it connected graded} if $A_0 = k$. 
Throughout this paper, the term $\mathbb{N}$-graded algebra refers to a locally finite $\mathbb{N}$-graded algebra.
By a graded $A$-module we mean a graded left $A$-module. For notation used throughout this paper but not introduced in this section, we refer the reader to \cite{WY}.

The {\it supremum degree} of a graded $A$-module (or, more generally, a graded vector space) $M$ is the supremum of the degrees of its nonzero homogeneous elements,
namely, 
$$\sdeg(M)=\sup \{ i \mid M_i \neq 0 \} \in \mathbb{Z} \cup \{ \pm \infty \}.$$ 
Similarly, the {\it infimum degree} of $M$ is the infimum of the degrees of its nonzero homogeneous elements: 
$$ \ideg(M)=\inf \{ i \mid M_i \neq 0 \} \in \mathbb{Z} \cup \{ \pm \infty \}.$$
By convention, $\sup(\emptyset)=-\infty$ and $\inf(\emptyset)=+\infty$.

Let $\xi=(\xi_0,\xi_1) \in \mathbb{R}\times \mathbb{R}$. For a graded vector space complex $X$ or for $X \in \mathbf{D}(\Gr{A})$, the {\it $\xi$-supremum degree} of $X$ is defined by
$$\sdeg_{\xi}(X)=\sup\{ \xi_1i + \xi_0 j \mid H^i(X)_j\neq 0,i,j\in \mathbb{Z} \}.$$

Similarly, for a graded vector space complex $X$ or for $X \in \mathbf{D}(\Gr{A})$, the {\it $\xi$-infimum degree} of $X$ is defined by
  $$\ideg_{\xi}(X)=\inf\{ \xi_1i+ \xi_0 j \mid H^i(X)_j\neq 0, i,j\in \mathbb{Z} \}.$$

Let $\D(-):=\gHom_k(-,k)$ denote the Matlis dual functor, and let $(-)^*:=\gHom_A(-,A)$.
The Matlis dual functor induces a duality between $\mathbf{D}_{lf}(\Gr A)$ and $\mathbf{D}_{lf}(\Gr A^o)$ (see, for example, \cite[Proposition 3.1]{VdB}).

For any $X \in \mathbf{D}(\Gr{A})$, we have
\begin{align}
\sdeg_{\xi}(X)=-\ideg_{\xi}(\D(X))=-\ideg_{-\xi}(X) \textrm{ and }\\
 \ideg_{\xi}(X)=-\sdeg_{\xi}(\D(X))=-\sdeg_{-\xi}(X).
\end{align}
We occasionally use the notation
\begin{align}
   \sup(X)=\mathop{\sup}\{ i\in \mathbb{Z} \mid  H^i(X)\neq 0 \} \textrm{ and } 
   \inf(X)=\mathop{\inf}\{  i\in \mathbb{Z}\mid  H^i(X)\neq 0 \}.
\end{align}

For any graded $A$-module $M$, the submodule
$\Gamma_A(M)= \mathop{\lim}\limits_{n\to \infty} \gHom_A(A/A_{\geqslant n},M)$ is called the {\it torsion submodule} of $M$.
The functor $\Gamma_A : \Gr A \to \Gr A$ is left exact. Let $R\Gamma_A$ denote its right derived functor. The $i$th {\it local cohomology} of $M$ is defined by $R^i\Gamma_A(M):=H^i(R\Gamma_A(M))\cong \mathop{\lim}\limits_{n\to \infty} \gExt^i_A(A/A_{\geqslant n},M)$. Similarly, we define $\Gamma_{A^o}$ and $R\Gamma_{A^o}$.

We say that $M\in \Gr A$  has {\it cohomological dimension} $d$ if $\sup(R\Gamma_A(M))=d$; we denote this dimension by $\cd(M)$. We say that the algebra $A$ or $\Gamma_A$ has {\it cohomological dimension} $d$ if $\sup\{ \cd(M)\mid M\in \Gr A\}=d;$ we denote this dimension by $\cd(\Gamma_A)$.

 \begin{definition} \label{Defi-CM}
    Suppose $M$ is a finitely generated graded left $A$-module. We call $M$ a graded Cohen-Macaulay module of dimension $s$ if $R^i\Gamma_A(M)=0$ for all $i\neq s$ and $R^s\Gamma_A(M)\neq 0$. We say that $A$ is a graded CM-algebra of dimension $s$ if ${}_{A}A$ is a graded Cohen-Macaulay module of dimension $s$.

    Obviously, $M$ is a graded Cohen-Macaulay module of dimension $s$ if and only if $\depth{M}=\cd(M)=s$.
\end{definition}

(Balanced) dualizing complexes were first defined in the noncommutative setting in \cite{Ye}. We refer the reader to \cite[Definition 2.12]{WY} for the definition. Van den Bergh proved an existence theorem for balanced dualizing complexes and a local duality theorem in \cite[Theorem 5.1, 6.3]{VdB}. These results also hold for locally finite graded algebras. If $R$ is a dualizing complex of $A$, then $R\gHom_A(-,R)$ and $R\gHom_{A^o}(-,R)$ define a duality between $\mathbf{D}_{fg}(\Gr A)$ and $\mathbf{D}_{fg}(\Gr A^o)$ (see \cite[Proposition 1.3]{YZ}).

$\mathbb{N}$-graded AS-Gorenstein algebras (sometimes called generalized AS-Gorenstein algebras) have been defined in several slightly different ways \cite{MV1,MV2,MM,RR,RRZ}. 
The following definition is given in \cite{LW}.
\begin{definition} \label{GAS-Gorenstein} 
   A noetherian $\mathbb{N}$-graded algebra $A$ is called an $\mathbb{N}$-graded Artin-Schelter Gorenstein algebra (or $\mathbb{N}$-graded AS-Gorenstein for short) of dimension $d$ if it satisfies the following conditions.
\begin{itemize}
    \item[(1)] $\injdim({}_{A}A)=\injdim(A_A)=d < \infty$,
    \item[(2)] For every graded simple left $A$-module $M$, $\gExt_A^i(M,A)=0$ if $i\neq d$; for every graded simple right $A$-module $N$, $\gExt_{A^o}^i(N,A)=0$ if $i\neq d$,
    \item[(3)] $\gExt_A^d(-,A)$ induces a bijection between the isomorphism classes of graded simple $A$-modules and those of graded simple $A^o$-modules, with inverse given by $\gExt_{A^o}^d(-,A)$.
\end{itemize}
Furthermore, if $A$ has finite global dimension, then $A$ is called an $\mathbb{N}$-graded Artin-Schelter regular algebra (or $\mathbb{N}$-graded AS-regular for short).
\end{definition}

\label{notation-e-i} Let $A$ be an $\mathbb{N}$-graded algebra and $S=A/J$ where $J$ is the graded Jacobson radical. Let $\{ S_1,S_2,\cdots,S_n\}$ be a set of representatives of all isomorphism classes of simple $S$-modules and, up to shifting, of all isomorphism classes of graded simple $A$-modules. There exists a set of orthogonal primitive idempotents $e_1,e_2,\cdots,e_n$ in $A_0$ such that $S_i\cong S\Bar{e_i}$, where $\Bar{e_i}$ is the image of $e_i$ in $A_0/J(A_0)$. If $e_1+e_2+\cdots+e_n=1$, then $A$ is called basic.

With the above notation, if $A$ is an $\mathbb{N}$-graded AS-Gorenstein algebra of dimension $d$, then there exists a permutation $\sigma \in \mathfrak{S}_n$ such that $\gExt^d_A(S_i,A)\cong S^{\prime}_{\sigma(i)}(\ell_i)$ for some $\ell_i\in \mathbb{Z}$ ($1\leqslant i\leqslant n$). Thus,
\begin{align}\label{def-GAS-Goren}
   \gExt^d_A(S,A)\cong \mathop{\bigoplus}\limits_{i=1}^n (e_{\sigma(i)}S(\ell_i))^{r_i}.
\end{align}
The set $\{ \ell_1,\ell_2,\cdots,\ell_n\}$ is called the set of Gorenstein parameters of $A$.
If $A$ is an $\mathbb{N}$-graded AS-Gorenstein algebra of dimension $d$, then $A$ is a CM-algebra of dimension $d$. 

\begin{definition}\cite[Definition 4.3]{IKU}
  $\ell_{av}^{A}:=n^{-1}\sum_{i=1}^n \ell_i \in \mathbb{Q}$ is called the average Gorenstein parameter of $A$.
\end{definition}

The following theorem summarizes several properties of $\mathbb{N}$-graded AS-Gorenstein algebras.
\begin{theorem}\cite[Lemma 3.8 and Theorem 4.12]{LW}\label{def-GAS}
Suppose $A$ is an $\mathbb{N}$-graded AS-Gorenstein algebra of dimension $d$. Using the notation above, $\gExt^d_{A^o}(S^{\prime}_i,A)\cong S_{\sigma^{-1}(i)}(\ell_{\sigma^{-1}(i)})$ for any $1\leqslant i\leqslant n$. Consequently,
\begin{align}\label{def-GAS-Goren coro}
   \gExt^d_{A^o}(S,A)\cong \mathop{\bigoplus}\limits_{i=1}^n (Se_{\sigma^{-1}(i)}(\ell_{\sigma^{-1}(i)}))^{r_i}.
\end{align}
Furthermore, $A$ admits a balanced dualizing complex given by $\D(R^d\Gamma_A(A))[d]$.

In fact, the balanced dualizing complex $R$ of $A$ has the form 
    \begin{align}   \label{GAS-Goren} 
     R=\D(R^d\Gamma_A(A))[d]\cong \mathop{\bigoplus}\limits_{i=1}^n (Ae_{\sigma(i)}(-\ell_i))^{r_i}[d]\cong \mathop{\bigoplus}\limits_{i=1}^n (e_iA(-\ell_i))^{r_{\sigma(i)}}[d].
     \end{align}
\end{theorem}

For $\mathbb{N}$-graded algebras, Minamoto and Mori also defined {\it AS-Gorenstein (AS-regular) property over $A_0$} \cite{MM}.  
\begin{definition} \cite[Definition 3.1]{MM}\label{def-AS-Goren over A_0}
    A noetherian $\mathbb{N}$-graded algebra $A$ is called AS-Gorenstein over $A_0$ (resp. AS-regular over $A_0$) of dimension $d$ with Gorenstein parameter $\ell$, if the following statements hold.
    \begin{itemize}
        \item [(1)] $\injdim{{}_AA}=d$ (resp. $\gldim A =d $);
        \item [(2)] there is an algebra automorphism $\nu$ over $A_0$ such that 
        $$R\gHom_A(A_0,A)\cong {}^1\D(A_0)^{\nu}(\ell)[-d]$$
        in $\mathbf{D}(\Gr A^e)$.
    \end{itemize}
\end{definition}

If $A$ is AS-Gorenstein of dimension $d$ over $A_0$, then $A$ is $\mathbb{N}$-graded AS-Gorenstein of dimension $d$. Furthermore, if $A$ is basic and $A_0$ is semisimple, then $A$ is AS-Gorenstein of dimension $d$ over $A_0$ if and only if $A$ is $\mathbb{N}$-graded AS-Gorenstein of dimension $d$.

\section{weighted homological regularities}\label{relationship between the regularities1}
Weighted regularities were introduced by Kirkman-Won-Zhang in \cite{KWZ2} for noetherian connected graded algebras. In this section, we extend their framework to locally finite $\mathbb{N}$-graded algebras and study weighted numerical regularities, including weighted Castelnuovo-Mumford regularities, weighted Tor-regularities (Ext-regularities), and weighted Ex-regularities. These include all the numerical regularities discussed in \cite{WY}. We also investigate relationships among these weighted numerical regularities. Homological invariants of complexes, including the supremum, infimum, projective dimension, depth, and $a$-invariant, fit naturally within the framework of weighted regularities for complexes.
 
The proofs of the results stated without proof in \S\ref{relationship between the regularities1} are analogous to those for the unweighted version in \cite{WY}, with the weight $(1,1)$ replaced by an arbitrary pair $\xi=(\xi_0,\xi_1)\in \mathbb{R}\times \mathbb{R}$. 

\subsection{Various numerical weighted regularities}
 In this subsection, we recall the definitions of the regularities mentioned above and explore their relationships in terms of the generating degrees of the minimal graded projective resolution of $X\in \mathbf{D}^-(\gr A)$, the Matlis duality functor $\D(-)$, and the duality functor $R\gHom_A(-,R)$ induced by the balanced dualizing complex $R$.

\begin{definition} Let $A$ be an $\mathbb{N}$-graded algebra , $S=A/J$ and $\xi=(\xi_0,\xi_1)\in \mathbb{R}\times \mathbb{R}$.
\begin{itemize}
   \item [(1)] \cite{KWZ2} \label{CM-reg-def} The $\xi$-Castelnuovo-Mumford regularity ({\it $\xi$-CM-regularity} for short) of $X\in \mathbf{D}^+(\Gr A)$ is defined to be
$$\CMreg_{\xi}(X)=\sdeg_{\xi} (R\Gamma_A(X)) =\sup_{i,j\in \mathbb{Z}} \{ \xi_1 i + \xi_0 j  \mid R^i\Gamma_A(X)_j\neq 0\}. $$
    \item [(2)] \cite{KWZ2} \label{torreg-def}  
 The {\it $\xi$-Tor-regularity} of $X\in \mathbf{D}^-(\Gr A)$ is defined to be 
$$\Torreg_{\xi}(X)=\sdeg_{\xi} (S \,\, {}^L\!\otimes_A X)=\sup_{i,j\in \mathbb{Z}}\{ -\xi_1 i + \xi_0 j  \mid \Tor^A_i(S,X)_j\neq 0\}.$$
\item [(3)] \cite{KWZ2} \label{Ext-reg-def}
  The {\it $\xi$-Ext-regularity} of $X\in \mathbf{D}^-(\Gr A)$ is defined to be $$\Extreg_{\xi}(X)=-\ideg_{\xi}(R\gHom_A(X, S))=-\inf_{i,j\in \mathbb{Z}} \{ \xi_1 i + \xi_0 j  \mid \gExt_A^i(X,S)_j\neq 0\}.$$ 
\item [(4)] \label{exreg-def}
         The {\it $\xi$-Ex-regularity} of $X\in \mathbf{D}^+(\Gr A)$ is defined to be 
    $$\Exreg_{\xi}(X)=\sdeg_{\xi} (R\gHom_A(S, X))=\sup_{i,j\in \mathbb{Z}}\{\xi_1i+\xi_0j \mid \gExt^i_A(S, X)_j\neq 0\}.$$ 
\end{itemize}

\end{definition}

\begin{remark}\label{relation-nonweighted-regularity}
 \begin{itemize}
   \item [(1)] In fact, for any $X\in \mathbf{D}^+(\Gr A)$, $\CMreg(X) = \CMreg_{(1,1)}(X)$, $\cmreg(X)=\CMreg_{(-1,-1)}(X)$, $\sup(R\Gamma_A(X))=\CMreg_{(0,1)}(X)$, and  \\
   $-\inf(R\Gamma_A(X))=-\depth{X}=\CMreg_{(0,-1)}(X)$, where $\cmreg(X)$ is defined in \cite[Definition 3.1]{WY}.
 \item [(2)] For any $X\in \mathbf{D}^-(\Gr A)$, $\Torreg(X)=\Torreg_{(1,1)}(X)$ and $\torreg(X)=\Torreg_{(-1,-1)}(X)$, where $\torreg(X)$ is defined in \cite[Definition 3.3]{WY}. 
 
 If $X\in \mathbf{D}^-(\Gr A)$ has a minimal projective resolution $P^{\bullet}$, then  $$\Torreg_{(0,1)}(X)=\sup(S \,\, {}^L\!\otimes_A X)=\sup(P^{\bullet}),$$  $$\Torreg_{(0,-1)}(X)=-\inf(S \,\, {}^L\!\otimes_A X)=\pdim{X}.$$
  \item [(3)]  For any $X\in \mathbf{D}^-(\Gr A)$, $\Extreg(X)=\Extreg_{(1,1)}(X)$ and  $\extreg(X)=\Extreg_{(-1,-1)}(X)$, where $\extreg(X)$ is defined in \cite[Definition 3.5]{WY}.
  
  If $X\in \mathbf{D}^-(\Gr A)$ has a minimal projective resolution $P^{\bullet}$, then $$\Extreg_{(0,1)}(X)=-\inf(R\gHom_A(X, S))=\sup(P^{\bullet}),$$ 
  $$\Extreg_{(0,-1)}(X)=\sup(R\gHom_A(X, S))=\pdim{X}.$$
 \item [(4)]  For any $X\in \mathbf{D}^+(\Gr A)$, $\exreg(X)=\Exreg_{(-1,-1)}(X)$, $\Exreg(X)=\Exreg_{(1,1)}(X)$, $\sup(R\gHom_A(S, X))=\Exreg_{(0,1)}(X)$ and 
 $$-\depth{X}=-\inf(R\gHom_A(S, X))=\Exreg_{(0,-1)}(X),$$ 
 where $\Exreg(X)$ and $\exreg(X)$ are defined in \cite[Definition 3.8]{WY}. 
 \end{itemize}
 
\end{remark}
\begin{remark}\label{a and b-invariant}
    \begin{itemize}
        \item [(1)] The $a$-invariant is given by
        $$a^*(X):= \sup_i \{ j \mid R^i\Gamma_A(X)_j \neq 0\}= \CMreg_{(1,0)}(X),$$ 
        where, in the commutative case, $a^*(X)$ is defined in \cite{GW, Tr2, Tr3}.
        \item [(2)] The $b$-invariant is given by
        $$b^*(X):=\sup_i \{b_i(X)\}=\sup_i \{j \mid \Tor_i^A(S,X)_j\neq 0\}=\Torreg_{(1,0)}(X),$$ 
      where $b_i(X):=\sup\{j\mid \Tor_i^A(S,X)_j\neq 0\}$ is defined in \cite{Tr3} for graded commutative algebras.
        \item [(3)] $\Torreg_{\xi}(S_A)=\Torreg_{\xi}({}_AS)$ for all $\xi=(\xi_0,\xi_1)\in \mathbb{R}\times \mathbb{R}$. This common value is denoted by $\Torreg_{\xi}(S)$ in the remainder of this paper. 
    \end{itemize}
\end{remark}

We now examine the conditions under which $\CMreg_{\xi}(X)$ is finite for $\xi=(\xi_0,\xi_1)\in \mathbb{R}\times \mathbb{R}$.
When $\xi_0\geqslant 0$, the following lemma is analogous to \cite[Lemma 3.2]{WY}.
\begin{lemma}\label{CM-reg-finite}
Suppose $A$ is a noetherian $\mathbb{N}$-graded algebra with a balanced dualizing complex. Then
\begin{itemize}
    \item [(1)]  $\CMreg_{\xi}(X) < \infty$ for all $0\neq X \in \mathbf{D}^b(\gr A)$ and $\xi=(\xi_0,\xi_1)\in \mathbb{R}_{\geqslant 0}\times \mathbb{R}$.
  \item[(2)] $\CMreg_{\xi}({}_{A}A)=\CMreg_{\xi}(A_A)$ for all $\xi=(\xi_0,\xi_1)\in \mathbb{R}\times \mathbb{R}$, which is denoted by $\CMreg_{\xi}(A)$.
\end{itemize}
\end{lemma}

The following proposition shows that finite weighted CM-regularity for various weights suffices to ensure that $A$ admits a balanced dualizing complex.
\begin{proposition}\label{bdc and CM}
Let $A$ be a noetherian $\mathbb{N}$-graded algebra. If, for all $0\neq M \in \Gr A$ and all $0\neq N \in \Gr A^o$, $\CMreg_{\xi}(M)<\infty$ and $\CMreg_{\xi}(N)<\infty$ for both $\xi=(1,0)$ and $\xi=(0,1)$, then $A$ has a balanced dualizing complex.
\end{proposition}
\begin{proof}
 Let $\xi=(0,1)$. Then $\CMreg_{(0,1)}(M)=\sup(R\Gamma_A(M))=\cd(M)<\infty$ for all $M\in \Gr A$. We claim that $\Gamma_A$ has finite cohomological dimension if $\cd(M)<\infty$ for all $M\in \Gr A$.

 Proof of claim. Let $\sup(R\Gamma_A(A))=d <\infty$. For any $M\in \gr A$, let $P^{\bullet}$ be a minimal graded projective resolution of $M$, $K^{-s}:=\ker d^{-s}$ for all $s\geqslant 0$ and $K^1:=M$. Then there is a short exact sequence
 \[ 0\to K^{-s}\to P^{-s}\to K^{-s+1}\to 0.\]
 It induces a long exact sequence
 \[ \cdots \to R^i\Gamma_A(P^{-s})\to R^i\Gamma_A(K^{-s+1})\to R^{i+1}\Gamma_A(K^{-s})\to R^{i+1}\Gamma_A(P^{-s})\to\cdots.\]
 Hence for any $i>d$, 
 \[R^i\Gamma_A(M)\cong R^{i+1}\Gamma_A(K^0)\cong R^{i+2}\Gamma_A(K^{-1})\cong \cdots\cong R^{i+s}\Gamma_A(K^{-s+1})\cong \cdots.\] 
 If there exists $i>d$ such that for any $s\geqslant0$, $R^{i+s}\Gamma_A(K^{-s+1})\neq 0$, then $\cd(K^{-s+1})\geqslant i+s>d$ and $\cd(K^{-s})\geqslant i+1+s$. Let $K=\mathop{\bigoplus}\limits_{s}K^{-s}\in \Gr A$. Then $\cd(K)=\infty$, which is a contradiction. Hence for any $i> d$, there exists $s\geqslant 0$ such that $R^{i+s}\Gamma_A(K^{-s+1})= 0$, that is,  for any $i>d$, $R^{i}\Gamma_A(M)= 0$. Hence $\cd(M)\leqslant d$ for all $M\in \gr A$. It follows from \cite[Lemma 4.3]{VdB} that $R^i\Gamma_A(-)$ commutes with direct limits. So for any $M\in \Gr A$, it follows from  $R^i\Gamma_A(M)=\varinjlim_j R^i\Gamma_A(M_j)$, where $M_j\in \gr A$ that $\cd(M)\leqslant d$. Hence $\cd(\Gamma_A)\leqslant d$.
 
 For $\xi=(1,0)$, we have $\CMreg_{(1,0)}(M)=\sup\{ j  \mid R^i\Gamma_A(M)_j\neq 0, i,j\in \mathbb{Z} \}<\infty$. Thus, $R^i\Gamma_A(M)$ is right bounded for any $M\in \gr A$ and any $i$. By \cite[Corollary 3.6 and Proposition 3.11]{AZ},
    $A$ satisfies the $\chi$-condition. 
    
    The argument for right modules is similar. Therefore, it follows from \cite[Theorem 6.3]{VdB} that $A$ has a balanced dualizing complex.
    \end{proof}
    
For $\xi_0\leqslant0$, the following proposition establishes that $\CMreg_{\xi}(X)< \infty$ if and only if the algebra $A$ is finite-dimensional. The result of \cite[Corollary 5.15]{WY} is recovered as a special case by taking $\xi=(-1,-1)$.

\begin{proposition}\label{finite-dim alg}
  Let $A$ be a noetherian $\mathbb{N}$-graded algebra with a balanced dualizing complex. The following conditions are equivalent:
   \begin{itemize}
       \item [(1)] $A$ is finite-dimensional.
       
       \item [(2)] $\CMreg_{\xi}(X)< \infty$ for all $0\neq X \in \mathbf{D}^b(\gr A)$ and all $\xi=(\xi_0,\xi_1)\in \mathbb{R}_{\leqslant 0}\times \mathbb{R}$.
       \item [(3)] $\CMreg_{\xi}(X)< \infty$ for all $0\neq X \in \mathbf{D}^b(\gr A)$ and some $\xi=(\xi_0,\xi_1)\in \mathbb{R}_{< 0}\times \mathbb{R}$.
   \end{itemize}  
\end{proposition}
\begin{proof}
(1)$\Rightarrow$(2), (2)$\Rightarrow$(3) Obviously.

 (3)$\Rightarrow$(1) For any $0 \neq X \in \mathbf{D}^b(\gr A)$, 
   we have
   \begin{align*}
       \CMreg_{\xi}(X)=&\sdeg_{\xi}(R\Gamma_A(X))\\
       =&-\ideg_{\xi}(\D(R\Gamma_A(X)))\\
       =&\sdeg_{-\xi} (\D(R\Gamma_A(X))) <\infty.
   \end{align*}
 
   This is equivalent to the condition
$\sdeg_{-\xi} (Y) < \infty$ for all $0 \neq Y \in \mathbf{D}^b(\gr A^o)$; hence, $A$ is finite-dimensional.
\end{proof}

For any $X\in \mathbf{D}^-(\Gr A)$ and $\xi=(\xi_0,\xi_1)\in \mathbb{R}\times \mathbb{R}$, the $\xi$-Tor-regularity and $\xi$-Ext-regularity of $X$ coincide, as stated in the following lemma (see \cite[Lemma 3.6]{WY}). Thus, we do not distinguish between $\Extreg_{\xi}(X)$ and $\Torreg_{\xi}(X)$ in this case.

\begin{lemma}\label{Extreg-Torreg}
Suppose that $A$ is an $\mathbb{N}$-graded algebra, $X\in \mathbf{D}^-(\Gr A)$, and $\xi=(\xi_0,\xi_1)\in \mathbb{R}\times \mathbb{R}$. Then $$\Extreg_{\xi}(X)=\Torreg_{\xi}(X).$$
\end{lemma}
 
There are more concrete descriptions of $\Torreg_{\xi}(X)\, (=\Extreg_{\xi}(X))$ by using the generating degrees of the minimal graded projective resolution of $0\neq X \in \mathbf{D}^-(\gr A)$. Suppose $P^{\bullet}$ is a minimal graded projective resolution of $X$. If $P^{-m}\neq 0$, we may assume that $P^{-m}$ has the following form (see \S\ref{notation-e-i})
$$P^{-m}= \mathop{\bigoplus}\limits_i \big(\mathop{\bigoplus}\limits_j Ae_i(-s_m^{i,j})\big).$$
We fix the notation that 
\begin{align}\label{u-l-degree-of-generators}
\uxm{m}{X}:=\sup \{ s_m^{i,j} \mid i, j\},\, \textrm{ and } \, \lxm{m}{X}:=\inf \{ s_m^{i,j} \mid i, j\}
\end{align}
which are the maximum degree and the minimal degree of the elements in the minimal generating set of $P^{-m}$ respectively. In fact, $\lxm{m}{X}=\ideg(P^{-m}).$

\begin{lemma}\label{notations of X}  Let $A$ be a noetherian $\mathbb{N}$-graded algebra and $0\neq X \in \mathbf{D}^-(\gr A)$. Then 
\begin{itemize}
    \item [(1)] $\Torreg_{\xi}(X)=\Extreg_{\xi}(X)=\sup\{\xi_0\uxm{m}{X}-\xi_1 m\mid P^{-m}\neq 0\}$ for all $\xi=(\xi_0,\xi_1)\in \mathbb{R}_{\geqslant 0}\times \mathbb{R}$.
    \item [(2)] $\Torreg_{\xi}(X)=\Extreg_{\xi}(X)=\sup \{\xi_0 \lxm{m}{X}-\xi_1 m \mid P^{-m}\neq 0\}$ for all $\xi=(\xi_0,\xi_1)\in \mathbb{R}_{\leqslant 0}\times \mathbb{R}$.   
 \end{itemize}

\end{lemma}

Suppose that $A$ has a balanced dualizing complex. Then $\D(R\Gamma_A(-))$ induces a duality between $\mathbf{D}_{fg}(\Gr A)$ and $\mathbf{D}_{fg}(\Gr A^o)$. Consequently, weighted Ex-regularity and weighted Ext-regularity are dual concepts when $A$ has a balanced dualizing complex. 
\begin{lemma}\label{dual-ext-ex} Let $A$ be a noetherian $\mathbb{N}$-graded algebra with a balanced dualizing complex $R$ and $\xi=(\xi_0,\xi_1)\in \mathbb{R}\times \mathbb{R}$. Then 
\begin{itemize}
 \item [(1)] for any $X\in \mathbf{D}^+_{fg}(\Gr A)$, 
 $\Exreg_{\xi}(X)=\Extreg_{-\xi}(\D(R\Gamma_A(X))) ,$
\item [(2)] for any $X\in \mathbf{D}^-_{fg}(\Gr A)$,
$\Extreg_{\xi}(X)=\Exreg_{-\xi}(\D(R\Gamma_A(X))).$
\end{itemize}

\end{lemma}

The Matlis dual establishes a duality between $\mathbf{D}_{lf}(\Gr A)$ and $\mathbf{D}_{lf}(\Gr A^o)$, thereby inducing a dual relationship between weighted Ex-regularity and weighted Ext-regularity (or weighted Tor-regularity).
\begin{lemma}\label{dual-ext-ex-by0Matlis} 
Let $A$ be an $\mathbb{N}$-graded algebra such that $S$ is finite-dimensional and $\xi=(\xi_0,\xi_1)\in \mathbb{R}\times \mathbb{R}$. Then
 \begin{itemize}
 \item [(1)] for any $X\in \mathbf{D}^+_{lf}(\Gr A)$,
\begin{align}\label{Ex-Ext and ex-ext-3}
  \Exreg_{\xi}(X)=\Extreg_{-\xi}(\D(X)), 
\end{align}
\item [(2)] for any $X\in \mathbf{D}^-_{lf}(\Gr A)$,
\begin{align}\label{Ex-Ext and ex-ext-4}
  \Extreg_{\xi}(X)=\Exreg_{-\xi}(\D(X)).
\end{align}
\end{itemize}
\end{lemma}

The following lemma describes the relationships among the numerical weighted regularities of the modules in a short exact sequence. It follows from the corresponding long exact sequences (see \cite[Lemma 3.7]{KWZ2}), where $$\reg_{\xi} \in \{\CMreg_{\xi}, \Torreg_{\xi}, \Extreg_{\xi},  \Exreg_{\xi}\}.$$

\begin{lemma}\label{exact sequence}
Let $A$ be an $\mathbb{N}$-graded algebra and $\xi=(\xi_0,\xi_1)\in \mathbb{R}_{\geqslant 0}\times \mathbb{R}$. For any distinguished triangle $X \to Y \to Z \to X[1]$ in $\mathbf{D}^b(\Gr A)$,
\begin{itemize}
    \item [(1)] $\reg_{\xi}(Y) \leqslant \max\{\reg_{\xi}(X),\reg_{\xi}(Z)\}$.
    \item [(2)] $\reg_{\xi}(X)\leqslant \max\{\reg_{\xi}(Y),\reg_{\xi}(Z)+\xi_1\}.$
    \item [(3)] $\reg_{\xi}(Z)\leqslant \max\{\reg_{\xi}(X)-\xi_1,\reg_{\xi}(Y)\}.$
\end{itemize}
\end{lemma}

\subsection{Further relationships between the regularities} \label{relationship between the regularities}
If $A_0$ is semisimple and $A$ has a balanced dualizing complex, we show in Theorem \ref{ex} that, under mild conditions,
$\CMreg_{\xi}(X) = \Exreg_{\xi}(X)$ for any $0\neq X \in \mathbf{D}^+_{fg}(\Gr A)$, generalizing \cite[Theorem 3.16]{WY} and \cite[Theorem 3.9]{Ngu} to the noncommutative graded case. We also prove in Proposition \ref{Extreg} that, if $A_0$ is semisimple and $\pdim{X}<\infty$, then $\Extreg_{\xi}(X)=-\ideg_{\xi}(R\gHom_A(X,A))$, generalizing \cite[Theorem 5.1]{Jo2} and \cite[Theorem 3.18]{WY}.
To prove Theorem \ref{ex} and Proposition \ref{Extreg}, we first prove Proposition \ref{CM and Ext}, which generalizes \cite[Lemma 5.2]{DW} and \cite[Proposition 3.14]{WY}. It states that if $A_0$ is semisimple, $\xi=(\xi_0,\xi_1)\in \mathbb{R}_{\leqslant 0}\times \mathbb{R}$, and $\xi_1\geqslant\xi_0$, then, for any $X\in \mathbf{D}^-(\gr A)$,
$$\Extreg_{\xi}(X)= \sdeg_{\xi}(X).$$

\begin{proposition}\label{CM and Ext}
Let $A$ be a noetherian $\mathbb{N}$-graded algebra, and let $\xi=(\xi_0,\xi_1)\in \mathbb{R}_{\leqslant 0}\times \mathbb{R}$. For any $0\neq X \in \mathbf{D}^-(\gr A),$
$$\Torreg_{\xi}(X)=\Extreg_{\xi}(X)\geqslant \sdeg_{\xi}(X).$$
Furthermore, if $A_0$ is semisimple and $\xi_1\geqslant\xi_0$, then 
$$\Torreg_{\xi}(X)=\sdeg_{\xi}(X).$$
\end{proposition}

\begin{remark}\label{remark of CM and Ext} 
   For any fixed integer $c\geqslant 0$, the argument in \cite[Proposition 3.14]{WY} shows that, for every $\xi=(\xi_0,\xi_1)\in \mathbb{R}_{\leqslant 0}\times \mathbb{R}$ satisfying $\xi_1\geqslant\xi_0$, 
   \begin{align*}
   &\sup \{ -\xi_1 i+ \xi_0 \ideg(H^{-i}(X)) \mid 0 \leqslant i\leqslant c\}\\
   =&-\inf \{\xi_1 i + \xi_0 \sdeg(\gExt_A^i(X,A_0)) \mid 0 \leqslant i\leqslant c\}.
   \end{align*}
This generalizes \cite[Corollary 5.2]{Jo2}, \cite[Theorem 5.2]{Sch}, and \cite[Remark 3.15]{WY}.
\end{remark}

Suppose that $A_0$ is semisimple and that $A$ has a balanced dualizing complex. The following theorem shows that weighted CM-regularity and weighted Ex-regularity coincide for any $0\neq X \in \mathbf{D}^+_{fg}(\Gr A)$, generalizing \cite[Theorem 3.9]{Ngu} and \cite[Theorem 3.16]{WY}. Proposition \ref{CM and Ext} plays a key role in the proof.

\begin{theorem}\label{ex}
Suppose that $A$ is a noetherian $\mathbb{N}$-graded algebra and $\xi=(\xi_0,\xi_1)\in \mathbb{R}_{\geqslant 0}\times \mathbb{R}$. 
Then, for any $0\neq X \in \mathbf{D}^+_{fg}(\Gr A)$, 
$$\CMreg_{\xi}(X)\leqslant \Exreg_{\xi}(X).$$
 Furthermore, if $\xi_1\leqslant\xi_0$, $A_0$ is semisimple, and $A$ has a balanced dualizing complex $R$, then 
 $$\CMreg_{\xi}(X)=\Exreg_{\xi}(X).$$
\end{theorem}
\begin{proof}
The proof of the inequality between $\CMreg_{\xi}(X)$ and $\Exreg_{\xi}(X)$ is analogous to that of \cite[Lemma 3.11]{WY}.

Let $Z=\D(R\Gamma_A(X))\in \mathbf{D}^-(\gr A^o)$.
Then, by Definition \ref{CM-reg-def},
$$\CMreg_{\xi}(X)=\sdeg_{\xi} (R\Gamma_A(X)) = -\ideg_{\xi}(Z)=\sdeg_{-\xi}(Z).$$
 
Since $R$ is a balanced dualizing complex of $A$, $R\gHom_A(A_0,X)\cong R\gHom_{A^o}(Z,A_0)$.
 By Definition \ref{exreg-def},
\begin{align*}
 \Exreg_{\xi}(X) &=\sdeg_{\xi} (R\gHom_A(A_0, X))\\
 &=\sdeg_{\xi} (R\gHom_{A^o}(Z, A_0))=\Extreg_{-\xi}(Z).   
\end{align*}
Since $-\xi=(-\xi_0,-\xi_1)$ satisfies $-\xi_0\leqslant 0$ and $-\xi_0\leqslant -\xi_1$, the right version of Proposition \ref{CM and Ext} implies that
\[\CMreg_{\xi}(X)=\sdeg_{-\xi}(Z)=\Extreg_{-\xi}(Z)=\Exreg_{\xi}(X). \qedhere \]
\end{proof}

By Lemma \ref{CM-reg-finite} and Theorem \ref{ex}, if $A$ has a balanced dualizing complex and $A_0$ is semisimple, then $\Exreg_{\xi}({}_{A}A)=\Exreg_{\xi}(A_A)$ for any $\xi=(\xi_0,\xi_1)\in \mathbb{R}_{\geqslant 0}\times \mathbb{R}$ with $\xi_1\leqslant\xi_0$. This common value is sometimes denoted by $\Exreg_{\xi}(A)$. If $A$ has a balanced dualizing complex $R$, then Lemma \ref{dual-ext-ex} (1) implies that $\Exreg_{-\xi}({}_{A}A)=\Torreg_{\xi}(R_A)$. Consequently, for any $\xi=(\xi_0,\xi_1)\in \mathbb{R}\times \mathbb{R}$, $\Exreg_{-\xi}({}_{A}A)=\Exreg_{-\xi}(A_A)$ if and only if $\Torreg_{\xi}({}_{A}R)=\Torreg_{\xi}(R_A)$.

The following proposition provides a weighted version of \cite[Theorem 3.18]{WY} and generalizes \cite[Theorem 5.1]{Jo2}.

\begin{proposition}\label{Extreg}
Suppose $A$ is a noetherian $\mathbb{N}$-graded algebra and $\xi=(\xi_0,\xi_1)\in \mathbb{R}_{\geqslant 0}\times \mathbb{R}$. For any $0\neq X\in \mathbf{D}^-(\gr A)$, $$\Torreg_{\xi}(X)=\Extreg_{\xi}(X) \geqslant -\ideg_{\xi}(R\gHom_A(X,A)).$$
 Furthermore, if $A_0$ is semisimple and $\xi_1\leqslant\xi_0$, then for any $0\neq X\in \mathbf{D}^b(\gr A)$ with $\pdim X < \infty$, $$\Torreg_{\xi}(X)=-\ideg_{\xi}(R\gHom_A(X,A)).$$
\end{proposition}

\section{Weighted numerical Artin-Schelter regularities}
\subsection{Weighted supremum and infimum degrees of Hom and Tensor complexes} \label{degree-of-Hom-Tensor}
In this subsection, we establish several inequalities and equalities concerning the weighted supremum or infimum degrees of $R\gHom_A(X,Y)$ and $Y\, {}^L\!\otimes_A X$ for $X, Y \in \mathbf{D}(\Gr A)$. Throughout this paper, saying that $A$ is a noetherian $\mathbb{N}$-graded algebra means that $A$ is a noetherian locally finite $\mathbb{N}$-graded algebra; hence, $\mathbf{D}(\gr A)$ is a full subcategory of $\mathbf{D}_{lf}(\Gr A)$.

The proof of the following proposition uses cohomological spectral sequences and follows from \cite[Proposition 4.1, 4.2]{WY}.

\begin{proposition} \label{sdeg-ideg inequality 1}
    Let $A$ be a noetherian $\mathbb{N}$-graded algebra and $\xi=(\xi_0,\xi_1)\in \mathbb{R}_{\geqslant 0}\times \mathbb{R}$. Then the following statements hold. 
    \begin{itemize}
        \item [(1)] For any $0\neq X\in \mathbf{D}^-(\gr A)$ and $0\neq Y\in \mathbf{D}^+(\Gr A)$, 
   \begin{align}  \label{inequality-1}
   -\ideg_{\xi}(R\gHom_A(X,Y))\leqslant \Extreg_{\xi}(X)-\ideg_{\xi}(Y).
   \end{align}
          \item [(2)] For any $0\neq X\in \mathbf{D}^-(\Gr A)$ with 
         torsion cohomologies and
        $0\neq Y\in \mathbf{D}^+(\Gr A)$,  
\begin{align}  \label{inequality-3-1}
   \sdeg_{\xi}(R\gHom_A(X,Y))\leqslant -\ideg_{\xi}(X)+\Exreg_{\xi}(Y).
   \end{align}
   \item [(3)] If $A$ has a balanced dualizing complex, then for any $0\neq X \in \mathbf{D}^-_{fg}(\Gr A)$
   and $0\neq Y \in \mathbf{D}^+_{fg}(\Gr A)$,  
\begin{align} \label{inequality-3-2}
   \sdeg_{\xi}(R\gHom_A(X,Y))\leqslant \CMreg_{-\xi}(X)+\Exreg_{\xi}(Y).
   \end{align}
    \end{itemize}
    \end{proposition}

The following proposition shows that the inequalities for the weighted‑supremum and weighted‑infimum degrees of $R\gHom_A(X,Y)$ and $Y\, {}^L\!\otimes_A X$ hold for all weight pairs in the plane. It also generalizes \cite[Theorems 3.3, 3.5]{KWZ2} and \cite[Propositions 4.4, 4.5]{WY}.

 \begin{proposition}\label{sdeg-ideg inequality 11} 
 Let $A$ be a noetherian $\mathbb{N}$-graded algebra and $\xi=(\xi_0,\xi_1)\in \mathbb{R}\times \mathbb{R}$. Then the following hold.      
 \begin{itemize}
   \item [(1)] For any $0\neq X\in \mathbf{D}^-(\Gr A)$
        and $0\neq Y\in \mathbf{D}^-(\gr A^o)$,
\begin{align}\label{inequality-2-1}
-\ideg_{\xi}(R\gHom_A(X,\D(Y)))\leqslant \sdeg_{\xi}(X)+\Extreg_{\xi}(Y).
\end{align}
\item [(2)] If $A$ has a balanced dualizing complex, then for any $0\neq X \in \mathbf{D}^-_{fg}(\Gr A)$ and $0\neq Y \in \mathbf{D}^-(\gr A^o)$ with torsion cohomologies,
\begin{align}  \label{inequality-2-2}
   -\ideg_{\xi}(R\gHom_A(X,\D(Y)))\leqslant \CMreg_{\xi}(X)+\Extreg_{\xi}(Y).
   \end{align}
   \item [(3)] For any $0\neq X\in \mathbf{D}^-(\gr A)$, $0\neq Y\in \mathbf{D}^-(\Gr A^{o})$,
    \begin{align}\label{inequality-2}
    -\ideg_{\xi}(Y\, {}^L\!\otimes_A X) \leqslant \Extreg_{-\xi}(X)-\ideg_{\xi}(Y).
    \end{align}
        \item [(4)] For any $0\neq X\in \mathbf{D}^-(\Gr A)$, $0\neq Y\in \mathbf{D}^-(\gr A^{o})$,
    \begin{align} \label{inequality-3}
    -\ideg_{\xi}(Y\, {}^L\!\otimes_A X) \leqslant -\ideg_{\xi}(X)+\Extreg_{-\xi}(Y).
    \end{align}
    
    \end{itemize}
 \end{proposition}  
 \begin{proof}
     Suppose that $\xi_0\geqslant 0$. The proof of this proposition follows the arguments in \cite[Proposition 4.4 and 4.5]{WY}.

     Suppose that $\xi_0 \leqslant 0$.

     For (1), it follows from $\D(Y\, {}^L\!\otimes_A X)\cong R\gHom_A(X,\D(Y))$ and \ref{inequality-3} holds for $\xi_0\geqslant 0$ that 
     \begin{align*}
         -\ideg_{\xi}(R\gHom_A(X,\D(Y)))=&\sdeg_{-\xi}(R\gHom_A(X,\D(Y)))\\
         =&\sdeg_{-\xi}(\D(Y\, {}^L\!\otimes_A X))\\
         =&-\ideg_{-\xi}(Y\, {}^L\!\otimes_A X)\\
         \leqslant& -\ideg_{-\xi}(X)+\Extreg_{\xi}(Y)\\
         =&\sdeg_{\xi}(X)+\Extreg_{\xi}(Y).
     \end{align*}

     For (2), it follows from \cite[Lemma 4.3]{WY} that 
     \begin{align*}
         R\gHom_A(X,\D(Y))\cong R\gHom_A(R\Gamma_A(X),\D(Y))\cong \D(Y\, {}^L\!\otimes_A R\Gamma_A(X)). 
     \end{align*}
     
Since \ref{inequality-3} holds for $\xi_0\geqslant 0$, 
      \begin{align*}
          -\ideg_{\xi}(R\gHom_A(X,\D(Y)))=&\sdeg_{-\xi}(R\gHom_A(X,\D(Y)))\\
          =&\sdeg_{-\xi}(\D(Y\, {}^L\!\otimes_A R\Gamma_A(X)))\\
          =&-\ideg_{-\xi}(Y\, {}^L\!\otimes_A R\Gamma_A(X))\\
          \leqslant & -\ideg_{-\xi}(R\Gamma_A(X))+\Extreg_{\xi}(Y)\\
          = &\CMreg_{\xi}(X)+\Extreg_{\xi}(Y).
      \end{align*}

      For (3), it follows from  $\D(Y\, {}^L\!\otimes_A X)\cong R\gHom_A(X,\D(Y))$ and \eqref{inequality-1} holds for $\xi_0\geqslant 0$ that
         \begin{align*}
         -\ideg_{\xi}(Y\, {}^L\!\otimes_A X)=&
             \sdeg_{-\xi}(Y\, {}^L\!\otimes_A X)\\
             =&-\ideg_{-\xi}(\D(Y\, {}^L\!\otimes_A X))\\
            =&-\ideg_{-\xi}(R\gHom_A(X,\D(Y))) \\
            \leqslant & \Extreg_{-\xi}(X)-\ideg_{-\xi}(\D(Y))\\
            =& \Extreg_{-\xi}(X)-\ideg_{\xi}(Y).
         \end{align*}
       
(4) Similar to the proof of (3) by using the right version of \eqref{inequality-1}, or 
\begin{align*}
     -\ideg_{\xi}(Y\, {}^L\!\otimes_A X)=&\sdeg_{-\xi}(Y\, {}^L\!\otimes_A X)\\
     =&-\ideg_{-\xi}(\D(Y\, {}^L\!\otimes_A X))\\
            =&-\ideg_{-\xi}(R\gHom_A(X,\D(Y))) \\
            \leqslant & \Extreg_{-\xi}(Y)+\sdeg_{-\xi}(X) \quad (\textrm{by \eqref{inequality-2-1}})\\
            =&\Extreg_{-\xi}(Y)-\ideg_{\xi}(X). \qedhere
         \end{align*}
 \end{proof}

Taking $Y=S$ in Proposition \ref{sdeg-ideg inequality 11} (2) yields $\Extreg_{\xi}(X)\leqslant \CMreg_{\xi}(X)+ \Extreg_{\xi}(S)$, generalizing \cite[Theorem 3.3]{KWZ2}. If $A$ has a balanced dualizing complex, setting $X=R\Gamma_A(A)$ in Proposition \ref{sdeg-ideg inequality 11} (1) yields $\CMreg_{\xi}(Y)\leqslant \Extreg_{\xi}(Y)+ \CMreg_{\xi}(A)$, generalizing \cite[Theorem 3.5]{KWZ2}.

 Next, we explore the conditions under which the inequalities \eqref{inequality-1}, \eqref{inequality-3-1},
 \eqref{inequality-2}, and 
 \eqref{inequality-3}
 in Propositions \ref{sdeg-ideg inequality 1} and  \ref{sdeg-ideg inequality 11} hold with equality. The equality cases are analogous to those for the other inequalities in the same propositions. In addition, they provide weighted generalizations of the results in \cite{WY}.
 
 We assume that $A_0$ is semisimple and first give two technical lemmas. The subsequent argument follows that for the non-weighted versions in \cite{WY}. 

\begin{lemma}\label{pre for equality 1}
Suppose $A$ is a noetherian $\mathbb{N}$-graded algebra with $A_0$ semisimple and $\xi_1\leqslant \xi_0$, where $\xi=(\xi_0,\xi_1)\in \mathbb{R}_{> 0}\times \mathbb{R}$. Let $0\neq X\in \mathbf{D}^-(\gr A)$, and let $P^{\bullet}$ be the minimal graded projective resolution of $X$.

If $p=-\ideg_{\xi}(X)$, then some $P^{-\alpha}$ has an indecomposable direct summand generated by an element of degree $\xi_0^{-1}(\xi_1\alpha-p)$ lying in $\ker d^{-\alpha}_{P^{\bullet}}$.
\end{lemma}

\begin{lemma}\label{pre for equality 2} Suppose $A$ is a noetherian $\mathbb{N}$-graded algebra with $A_0$ semisimple and $\xi_1\leqslant \xi_0$, where $\xi=(\xi_0,\xi_1)\in \mathbb{R}_{> 0}\times \mathbb{R}$. Let $0\neq X\in \mathbf{D}^-(\gr A)$, and let $P^{\bullet}$ be the minimal graded projective resolution of $X$.

If $p=-\ideg_{\xi}(X)=-\ideg_{\xi}(e_kX)$ for some $1\leqslant k \leqslant n$, then there exists $\beta \in \mathbb{Z}$ with $H^{-\beta}(e_kP^{\bullet})\neq 0$ such that 
\[\ideg(\ker d_{P^{\bullet}}^{-\beta})=\xi_0^{-1}(\xi_1\beta-p)\, \textrm{ and } \, e_k(\ker d_{P^{\bullet}}^{-\beta})_{\xi_0^{-1}(\xi_1\beta-p)} \neq 0.\]
\end{lemma}

We now give criteria under which \eqref{inequality-1}, \eqref{inequality-3-1}, \eqref{inequality-2} and \eqref{inequality-3} hold with equality when $A_0$ is semisimple. Note that if $A_0$ is semisimple, $\xi=(\xi_0,\xi_1)\in \mathbb{R}_{\leqslant 0}\times \mathbb{R}$ and $\xi_1\geqslant\xi_0$, then $\sdeg_{\xi}(X) = \Extreg_{\xi}(X)$ by Proposition \ref{CM and Ext}. If $A$ is connected graded, then condition (2) in each of the following propositions and corollaries in this subsection is automatically satisfied, so the identities in (1) always hold in the connected graded case.

Propositions \ref{sdeg-ideg equality 1} and \ref{sdeg-ideg equality 2} provide equivalent conditions for equality in \eqref{inequality-2} and \eqref{inequality-3} when $A_0$ is semisimple.

\begin{proposition}\label{sdeg-ideg equality 1}
   Suppose $A$ is a noetherian $\mathbb{N}$-graded algebra with $A_0$ semisimple and $\xi_1\leqslant \xi_0$, where $\xi=(\xi_0,\xi_1)\in \mathbb{R}_{> 0}\times \mathbb{R}$. For any $0\neq Y \in \mathbf{D}^b(\gr {A^o})$, the following are equivalent.
 \begin{itemize}
     \item [(1)] $\ideg_{\xi}(Y\, {}^L\!\otimes_A X) =  \ideg_{\xi}(X)+\ideg_{\xi}(Y)$ for all $0\neq X\in \mathbf{D}^b(\gr A)$.
     \item [(2)] $\ideg_{\xi}(Ye_i)=\ideg_{\xi}(Y)$ for all $1\leqslant i\leqslant n$.
 \end{itemize} 
\end{proposition}
\begin{proof}
    (1) $\Rightarrow$ (2) It follows by taking $X=Ae_i$.

    (2) $\Rightarrow$ (1)
Without loss of generality, let $X^i=0, Y^i=0$ for all $i>0$ and $\sup(X)=\sup(Y)=0$. 
Since $0\neq X\in \mathbf{D}^b(\gr A)$ and $0 \neq Y \in \mathbf{D}^b(\gr {A^o})$, $-\ideg_{\xi}(X) \neq  \pm \infty$ and $-\ideg_{\xi}(Y) \neq \pm \infty$. Let $p_1= -\ideg_{\xi}(X)$ and $p_2=-\ideg_{\xi}(Y)$. 
Let $P^{\bullet}$ and $Q^{\bullet}$ be the minimal graded projective resolutions of $X$ and $Y$ respectively.

Note that $Y\, {}^L\!\otimes_A X\cong Q^{\bullet} \otimes_{A} P^{\bullet}$ 
    where $(Q^{\bullet} \otimes_{A} P^{\bullet})^m= \mathop{\bigoplus}\limits_{q}Q^{q} \otimes_{A} P^{m-q},$ and its $m$th-differential is the morphism $$\partial^m(x \otimes y)=d_{Q^{\bullet}}(x) \otimes y +(-1)^{q}x\otimes d_{P^{\bullet}}(y).$$
It follows from Proposition \ref{sdeg-ideg inequality 11} (3) that 
     \begin{align*}
         &-\ideg_{\xi}(Y\, {}^L\!\otimes_A X)\\
         &=-\inf\{\xi_1 m+\xi_0 \ideg(H^m(Q^{\bullet}\otimes_{A} P^{\bullet}))\mid \Tor_{-m}^A(Y,X)\neq 0\}\\
         &\leqslant \Extreg_{-\xi}(X)-\ideg_{\xi}(Y)\\
         &= \sdeg_{-\xi}(X) +p_2\\
         &=-\ideg_{\xi}(X)+p_2\\
         &=p_1+p_2.
     \end{align*}
Hence, for any $m\in \mathbb{Z}$ such that $\Tor_{-m}^A(Y,X)\neq 0$, $$\ideg(H^{m}(Q^{\bullet}\otimes_{A} P^{\bullet})) \geqslant \xi_0^{-1}(-p_1-p_2-\xi_1 m).$$
To prove (1), it suffices to show that there is some $\alpha \in \mathbb{Z}$ with
$\Tor_{-\alpha}^A(Y,X)\neq 0$ such that
$$\ideg(H^{\alpha}(Q^{\bullet}\otimes_{A} P^{\bullet}))\leqslant \xi_0^{-1}(-p_1-p_2-\xi_1 \alpha).$$

On one hand, it follows from Lemma \ref{pre for equality 1} that there is some $\alpha_1 \in \mathbb{Z}$ such that $P^{-{\alpha}_1}$ has an indecomposable direct summand generated by an element of degree $\xi_0^{-1}(\xi_1 \alpha_1-p_1)$ contained in $\ker d^{-\alpha_1}_{P^{\bullet}}$.  
So, there is some
$k$ with $1\leqslant k\leqslant n$ such that $Ae_k(\xi_0^{-1}(p_1-\xi_1\alpha_1))\subseteq \ker d^{-\alpha_1}_{P^{\bullet}}$.

On the other hand, by the assumption in (2) and the right version of Lemma \ref{pre for equality 2}, there exists $\alpha_2 \in \mathbb{Z}$ with $H^{-\alpha_2}(Q^{\bullet}e_k)\neq 0$ such that 
$$\ideg(\ker d_{Q^{\bullet}}^{-{\alpha_2}})=\xi_0^{-1}(\xi_1\alpha_2-p_2)\, \textrm{ and } \, (\ker d_{Q^{\bullet}}^{-{\alpha_2}})_{\xi_0^{-1}(\xi_1\alpha_2-p_2)}\,e_k \neq 0.$$
There is some $0\neq y \in (\ker d_{Q^{\bullet}}^{-\alpha_2})_{\xi_0^{-1}(\xi_1\alpha_2-p_2)} \subseteq Q^{-\alpha_2}$ such that 
$$ 0 \neq y \otimes_A e_k(\xi_0^{-1}(p_1-\xi_1\alpha_1)) \subseteq Q^{-\alpha_2} \otimes_A P^{-\alpha_1}.$$

Then $0\neq y \otimes_A e_k(\xi_0^{-1}(p_1-\xi_1\alpha_1)) \in \ker \partial^{\alpha}$, where $\alpha=-\alpha_2-\alpha_1$. Hence $$\ideg(\ker \partial^{\alpha})\leqslant  \xi_0^{-1}(\xi_1\alpha_2-p_2) + \xi_0^{-1}(\xi_1\alpha_1 -p_1) = -\xi_1\alpha-p_1-p_2.$$ 
Then $0\neq y \otimes_A e_k(\xi_0^{-1}(p_1-\xi_1\alpha_1)) \in \ker \partial^{\alpha}$, where $\alpha=-\alpha_2-\alpha_1$. Hence $$\ideg(\ker \partial^{\alpha})\leqslant  \xi_0^{-1}(\xi_1\alpha_2-p_2) + \xi_0^{-1}(\xi_1\alpha_1 -p_1) = \xi_0^{-1}(-\xi_1\alpha-p_1-p_2).$$ 
Since $P^{\bullet}$ and $Q^{\bullet}$ are minimal graded projective resolutions, 
$$\im \partial^{\alpha-1} \subseteq \mathop{\bigoplus}\limits_{q}(Q^{q+1}A_{\geqslant 1} \otimes_{A} P^{\alpha-1-q}+Q^{q} \otimes_{A} A_{\geqslant 1}P^{\alpha-q}) \subseteq (Q^{\bullet}\otimes_{A} P^{\bullet})^{\alpha}.$$ 
It follows that 
\begin{align*}
&\ideg(\im \partial^{\alpha-1}) \\
\geqslant& \min\{\lxm{-q-1}{Y}+1 + \lxm{-\alpha+1+q}{X}, \lxm{-q}{Y}+1 + \lxm{-\alpha+q}{X}  \}\\
\geqslant& \min\{\xi_0^{-1}(\xi_1(-q-1)-p_2)+1 +\xi_0^{-1}(\xi_1(-\alpha+1+q)-p_1),\\ &\xi_0^{-1}(-\xi_1q-p_2)+1 +\xi_0^{-1}(\xi_1(-\alpha+q)-p_1) \}\\
=& \xi_0^{-1}(-p_1-p_2-\xi_1\alpha)+1.
\end{align*}
Hence
    $\ideg (H^{\alpha}(Q^{\bullet}\otimes_{A} P^{\bullet}))\leqslant \xi_0^{-1}(-p_1-p_2-\xi_1\alpha).$
    
Therefore, 
$-\ideg_{\xi}(Y\, {}^L\!\otimes_A X)=p_1+p_2=-\ideg_{\xi}(X)-\ideg_{\xi}(Y),$ and the proof of (2) $\Rightarrow$ (1) is finished.
\end{proof}

\begin{proposition}\label{sdeg-ideg equality 2}
   Suppose that $A$ is a noetherian $\mathbb{N}$-graded algebra with $A_0$ semisimple and  $\xi_1\leqslant \xi_0$, where $\xi=(\xi_0,\xi_1)\in \mathbb{R}_{> 0}\times \mathbb{R}$. 
   For any $0\neq X\in \mathbf{D}^b(\gr A)$, the following are equivalent.
 \begin{itemize}
     \item [(1)] $\ideg_{\xi}(Y\, {}^L\!\otimes_A X) =  \ideg_{\xi}(X)+\ideg_{\xi}(Y)$ for all $0\neq Y \in \mathbf{D}^b(\gr {A^o})$.
     \item [(2)] $\ideg_{\xi}(e_iX)=\ideg_{\xi}(X)$ for all $1\leqslant i\leqslant n$.
 \end{itemize} 
\end{proposition}

\begin{proof} (1) $\Rightarrow$ (2) It follows by taking $Y=e_iA$  in (1).

(2) $\Rightarrow$ (1) Similar to the proof of (2) $\Rightarrow$ (1) in Proposition \ref{sdeg-ideg equality 1}.
\end{proof}

Propositions \ref{sdeg-ideg equality 3} and \ref{sdeg-ideg equality 4} provide criteria for equality to hold in \eqref{inequality-1} when $A_0$ is semisimple.
\begin{proposition}\label{sdeg-ideg equality 3}
 Suppose that $A$ is a noetherian $\mathbb{N}$-graded algebra with $A_0$ semisimple and  $\xi_1\leqslant \xi_0$, where $\xi=(\xi_0,\xi_1)\in \mathbb{R}_{> 0}\times \mathbb{R}$. For any $0\neq Y \in \mathbf{D}^b(\gr {A})$, the following conditions are equivalent.
 \begin{itemize}
    \item [(1)] $-\ideg_{\xi}(R\gHom_A(X,Y))=\Extreg_{\xi}(X)-\ideg_{\xi}(Y)$ for all $0\neq X\in \mathbf{D}^b(\gr A)$ with $\pdim(X)<\infty$.
    \item [(2)] 
    $\ideg_{\xi}(e_iY)=\ideg_{\xi}(Y)$ for all $1\leqslant i\leqslant n$.
 \end{itemize}
\end{proposition} 
\begin{proof}
    (1) $\Rightarrow$ (2) It follows by taking $X=Ae_i$  in (1) that
    $$-\ideg_{\xi}(R\gHom_A(Ae_i,Y))=\Extreg_{\xi}(Ae_i)-\ideg_{\xi}(Y).$$
    Since $R\gHom_A(Ae_i,Y)\cong e_iY$ as graded $k$-spaces and $\Extreg_{\xi}(Ae_i)=0$, $$\ideg_{\xi}(e_iY)=0+\ideg_{\xi}(Y)=\ideg_{\xi}(Y).$$ 

    (2) $\Rightarrow$ (1) Let $Z=R\gHom_A(X,A)$. Since $\pdim(X)<\infty$, $0\neq Z \in \mathbf{D}^b(\gr A^o)$ and 
      \begin{align}\label{isomorphism}
          R\gHom_A(X,Y) \cong R\gHom_A(X,A) \, {}^L\!\otimes_A Y \cong Z \, {}^L\!\otimes_A Y. 
      \end{align}
      It follows from Proposition \ref{sdeg-ideg equality 2} that 
      $$ \ideg_{\xi}(Z \, {}^L\!\otimes_A Y)=\ideg_{\xi}(Y)+\ideg_{\xi}(Z).$$
      By Theorem \ref{Extreg}, 
      $\Extreg_{\xi}(X)=-\ideg_{\xi}(R\gHom_A(X,A))=-\ideg_{\xi}(Z).$

      Thus, 
      \begin{align*}
          -\ideg_{\xi}(R\gHom_A(X,Y))=&-\ideg_{\xi}(Z \, {}^L\!\otimes_A Y)\\
          =&-\ideg_{\xi}(Y)-\ideg_{\xi}(Z)\\
          =&-\ideg_{\xi}(Y)+\Extreg_{\xi}(X).
      \end{align*}
The proof of (2) $\Rightarrow$ (1) is finished.
\end{proof}

\begin{proposition}\label{sdeg-ideg equality 4}
 Suppose that $A$ is a noetherian $\mathbb{N}$-graded algebra with $A_0$ semisimple and $\xi_1\leqslant \xi_0$, where $\xi=(\xi_0,\xi_1)\in \mathbb{R}_{> 0}\times \mathbb{R}$. For any $0\neq X \in \mathbf{D}^b(\gr A)$ with $\pdim(X)<\infty$, the following are equivalent.
 \begin{itemize}
     \item [(1)] $-\ideg_{\xi}(R\gHom_A(X,Y))=\Extreg_{\xi}(X)-\ideg_{\xi}(Y)$ for all $0\neq Y\in \mathbf{D}^b(\gr A)$.
    \item [(2)] $\Extreg_{\xi}(X)=-\ideg_{\xi}(R\gHom_A(X,Ae_i))$ for all $1\leqslant i\leqslant n$. 
 \end{itemize}
\end{proposition} 

\begin{proof} (1) $\Rightarrow$ (2) It follows by taking $Y=Ae_i$.

(2) $\Rightarrow$ (1) Let $Z=R\gHom_A(X,A)$. Since $\pdim(X)<\infty$, $0\neq Z \in \mathbf{D}^b(\gr A^o)$. It follows from Theorem \ref{Extreg} and $R\gHom_A(X,Ae_i)\cong Z \otimes_A Ae_i\cong Ze_i$ that 
$$ -\ideg_{\xi}(Z)=\Extreg_{\xi}(X)\stackrel{(2)}{=}-\ideg_{\xi}(R\gHom_A(X,Ae_i))=-\ideg_{\xi}(Ze_i).$$ 
 Hence by Proposition \ref{sdeg-ideg equality 1} and \eqref{isomorphism}, 
 \begin{align*}
      -\ideg_{\xi}(R\gHom_A(X,Y))=&-\ideg_{\xi}(Z \, {}^L\!\otimes_A Y)\\
      =&-\ideg_{\xi}(Y)-\ideg_{\xi}(Z)\\
      =&-\ideg_{\xi}(Y)+\Extreg_{\xi}(X).
 \end{align*}
 The proof of (2) $\Rightarrow$ (1) is finished.
\end{proof}

Proposition \ref{sdeg-ideg equality 8} specifies when equality holds in \eqref{inequality-3-1}.

\begin{proposition}\label{sdeg-ideg equality 8}
    Suppose that $A$ is a noetherian $\mathbb{N}$-graded algebra with $A_0$ semisimple and  $\xi_1\leqslant \xi_0$, where $\xi=(\xi_0,\xi_1)\in \mathbb{R}_{> 0}\times \mathbb{R}$.  For any $0\neq X \in \mathbf{D}^b(\gr A)$, the following are equivalent.
 \begin{itemize}
      \item [(1)] $\sdeg_{\xi}(R\gHom_A(X,Y)) = -\ideg_{\xi}(X)+\Exreg_{\xi}(Y)$ for all $0\neq Y \in \mathbf{D}^b(\Gr {A})$ with $\D(Y) \in \mathbf{D}^b(\gr {A^o})$.
     \item [(2)] $\ideg_{\xi}(e_iX)=\ideg_{\xi}(X)$ for all $1 \leqslant i\leqslant n$.
 \end{itemize} 
\end{proposition}
 \begin{proof}
     (1) $\Rightarrow$ (2)  By  taking $Y=\D(e_iA)$, it follows from $\Exreg_{\xi}(Y)=\Extreg_{-\xi}(\D(Y))=0$.

     (2) $\Rightarrow$ (1) Since $\D(Y) \in \mathbf{D}^b(\gr {A^o})$, $Y \in \mathbf{D}^b_{lf}(\Gr {A})$ and $\D(\D(Y)\cong Y$.
      Then $\D(\D(Y)\, {}^L\!\otimes_A X)\cong R\gHom_A(X,\D(\D(Y))) \cong R\gHom_A(X,Y)$. 
     Hence
     \begin{align*}
         \sdeg_{\xi}(R\gHom_A(X,Y))=&-\ideg_{\xi}(\D(Y)\, {}^L\!\otimes_A X)\\
         =&-\ideg_{\xi}(X)-\ideg_{\xi}(\D(Y))  \quad (\textrm{by Proposition } \ref{sdeg-ideg equality 2})\\
         =& -\ideg_{\xi}(X)+\sdeg_{-\xi}(\D(Y))\\
         =&-\ideg_{\xi}(X)+\Extreg_{-\xi}(\D(Y))\\
        \stackrel{(\ref{Ex-Ext and ex-ext-3})}{=}&-\ideg_{\xi}(X)+\Exreg_{\xi}(Y). \qedhere
     \end{align*}
\end{proof}

\subsection{Weighted numerical Artin-Schelter regularity}
In this subsection, following \cite{KWZ2}, we first introduce the concept of weighted numerical Artin–Schelter regularity for noetherian 
$\mathbb{N}$-graded algebras. We then extend Theorems 3.3 and 3.5 of \cite{KWZ2} to 
$\mathbb{N}$-graded algebras equipped with balanced dualizing complexes and establish further analogous, and in some cases stronger, relationships among these regularities (see, for instance, Theorem \ref{relation between CM and Tor} (4) and Proposition \ref{relation between CM and Tor-p}).
For any noetherian $\mathbb{N}$-graded algebra $A$ with a balanced dualizing complex and any $\xi=(\xi_0,\xi_1)\in \mathbb{R}\times \mathbb{R}$, we prove that the following conditions are equivalent: $\ASreg_{\xi}(A)=0$; $\CMreg_{\xi}(X)= \Torreg_{\xi}(X)+\CMreg_{\xi}(A)$ for all $0\neq X \in \mathbf{D}^b(\gr A)$; and $\Torreg_{\xi}(X)= \CMreg_{\xi}(X)+\Torreg_{\xi}(S)$ for all $0\neq X \in \mathbf{D}^b(\gr A)$.

\begin{definition} (see \cite[Definition 0.2]{KWZ2})
Let $\xi\in \mathbb{R}\times \mathbb{R}$. For a noetherian $\mathbb{N}$-graded algebra $A$, the $\xi$-Artin-Schelter regularity ({\it$\xi$-AS regularity} for short) is defined by $$\ASreg_{\xi}(A) =\CMreg_{\xi}(A) + \Torreg_{\xi}(S).$$
\end{definition}

Obviously, $\ASreg_{(1,1)}(A)=\ASreg(A)$ and $\ASreg_{(-1,-1)}(A)=\asreg(A)$, where $\ASreg(A)$ and $\asreg(A)$ are defined in \cite[Definition 5.1 and 5.3]{WY}.

Suppose that $A$ is a noetherian $\mathbb{N}$-graded algebra with $A_0$ semisimple and that $A$ has a balanced dualizing complex. If $\xi_1\leqslant\xi_0$, where $\xi=(\xi_0,\xi_1)\in \mathbb{R}_{\geqslant 0}\times \mathbb{R}$, then $\ASreg_{\xi}(A)=\Exreg_{\xi}(A)+\Torreg_{\xi}(S)$ by Theorem \ref{ex}.

Parts (1) and (2) of Theorem \ref{relation between CM and Tor} generalize \cite[Theorems 3.3 and 3.5]{KWZ2}. Moreover, the condition on $\xi$ can be relaxed to allow $\xi$ to vary over the entire plane $\mathbb{R}\times \mathbb{R}$. The proofs rely on Proposition \ref{sdeg-ideg inequality 1} (1), \eqref{inequality-2}, and Proposition \ref{sdeg-ideg inequality 11} (2).

\begin{theorem}\label{relation between CM and Tor}
 Suppose that $A$ is a noetherian $\mathbb{N}$-graded algebra with a balanced dualizing complex $R$. Let $\xi=(\xi_0,\xi_1)\in \mathbb{R}\times \mathbb{R}$. Then
   \begin{itemize}
       \item [(1)] $\CMreg_{\xi}(X)\leqslant \Torreg_{\xi}(X)+\CMreg_{\xi}(A)$ for all $0\neq X \in \mathbf{D}^-(\gr A)$.
       \item [(2)] $\Torreg_{\xi}(X)\leqslant \CMreg_{\xi}(X)+\Torreg_{\xi}(S)$ for all $0\neq X \in \mathbf{D}^-(\gr A)$.
       \item [(3)] $\ASreg_{\xi}(A)\geqslant 0$.
       
   \end{itemize}
\end{theorem}

For $\xi=(1,0)$, Theorem \ref{relation between CM and Tor} yields several inequalities relating the $a$-invariant and the $b$-invariant, stated in the following corollary. The definitions of these invariants are given in Remark \ref{a and b-invariant}.

\begin{corollary}\label{relation-a-b}
    Suppose that $A$ is a noetherian $\mathbb{N}$-graded algebra with a balanced dualizing complex $R$. Then
   \begin{itemize}
       \item [(1)] $a^*(X)\leqslant b^*(X)+a^*(A)$ for all $0\neq X \in \mathbf{D}^-(\gr A)$.
       \item [(2)] $b^*(X)\leqslant a^*(X)+b^*(S)$ for all $0\neq X \in \mathbf{D}^-(\gr A)$.
       \item [(3)] $a^*(A)+b^*(S)\geqslant 0$.
   \end{itemize}
\end{corollary}

In the non-weighted or connected graded case, the condition that the inequality in (1) or (2) of Theorem \ref{relation between CM and Tor} is always an equality is closely related to the Artin-Schelter regularity or finite-dimensionality of the algebra (\cite[Corollary 3.4]{KWZ2} and \cite[Corollary 5.16]{WY}).

The following corollary shows that $\ASreg_{\xi}(A)=0$ characterizes when the inequality in (1) or (2) of Theorem \ref{relation between CM and Tor} is always an equality. For a proof, see \cite[Corollary 5.6 and 5.17]{WY}.

\begin{corollary}\label{ASreg=0}
    Suppose that $A$ is a noetherian $\mathbb{N}$-graded algebra with a balanced dualizing complex. Let $\xi=(\xi_0,\xi_1)\in \mathbb{R}\times \mathbb{R}$. Then the following conditions are equivalent.
    \begin{itemize}
        \item [(1)] $\ASreg_{\xi}(A)=0$.
        
        \item [(2)] $\CMreg_{\xi}(X)=\Torreg_{\xi}(X)+\CMreg_{\xi}(A)$ for all $0\neq X \in \mathbf{D}^b(\gr A)$.
        \item [(3)] $\Torreg_{\xi}(X)=\CMreg_{\xi}(X)+\Torreg_{\xi}(S)$ for all $0\neq X \in \mathbf{D}^b(\gr A)$.
        \item [(4)] There exists some $0\neq X \in \mathbf{D}^b(\gr A)$ such that $$\CMreg_{\xi}(X)=\Torreg_{\xi}(X)+\CMreg_{\xi}(A) \, \textrm { and } \, \Torreg_{\xi}(X)=\CMreg_{\xi}(X)+\Torreg_{\xi}(S).$$
    \end{itemize}
\end{corollary}

The following corollary establishes several new relationships between the $a$-invariant and the $b$-invariant by setting $\xi=(1,0)$ in Corollary \ref{ASreg=0}.

 \begin{corollary}\label{a-invariant+b-invariant=0}
    Suppose that $A$ is a noetherian $\mathbb{N}$-graded algebra with a balanced dualizing complex. Then the following are equivalent.
    \begin{itemize}
        \item [(1)] $a^*(A)+b^*(S)=0$.
        
        \item [(2)] $a^*(X)=b^*(X)+a^*(A)$ for all $0\neq X \in \mathbf{D}^b(\gr A)$.
        \item [(3)] $b^*(X)=a^*(X)+b^*(S)$ for all $0\neq X \in \mathbf{D}^b(\gr A)$.
        \item [(4)] There exists $0\neq X \in \mathbf{D}^b(\gr A)$ such that $$a^*(X)=b^*(X)+a^*(A) \, \textrm { and } \, b^*(X)=a^*(X)+b^*(S).$$
    \end{itemize}
\end{corollary}

The following corollary shows that weighted CM-regularity or weighted AS-regularity serves as an invariant measuring how far the algebra $A$ is from being finite-dimensional. 

\begin{corollary}\label{cmreg-finite}
 Suppose that $A$ is a noetherian $\mathbb{N}$-graded algebra with a balanced dualizing complex $R$. Let $\xi=(\xi_0,\xi_1)\in \mathbb{R}_{< 0}\times \mathbb{R}$ with $\xi_1\geqslant \xi_0$. If $A_0$ is semisimple, then the following conditions are equivalent. 
\begin{itemize}
        \item [(1)] $\CMreg_{\xi}(A)<\infty$.
        \item [(2)] $\CMreg_{\xi}(X)< \infty$ for all $0\neq X \in \mathbf{D}^b(\gr A)$.
        \item [(3)] $A$ is finite-dimensional.
        \item [(4)] $\CMreg_{\xi}(A)=0$.
        \item [(5)] $\ASreg_{\xi}(A)=0$.
        \item [(6)] $\ASreg_{\xi}(A)< \infty$.
\end{itemize}     
\end{corollary}

\begin{proof}  (1) $\Leftrightarrow$ (2) Since $A_0$ is semisimple, $\Extreg_{\xi}(S)=0$. 
Thus, $\CMreg_{\xi}(A)\geqslant 0$ by Theorem \ref{relation between CM and Tor} (3). By Proposition \ref{relation between CM and Tor} (1) and Proposition \ref{CM and Ext}, for any $0\neq X \in \mathbf{D}^b(\gr A)$,
$$\CMreg_{\xi}(X)\leqslant\Extreg_{\xi}(X)+\CMreg_{\xi}(A)=\sdeg_{\xi}(X)+\CMreg_{\xi}(A).$$
Since $\sdeg_{\xi}(X)<\infty$, it follows that $\CMreg_{\xi}(A)<\infty$ if and only if $\CMreg_{\xi}(X)< \infty$ for all $0 \neq X \in \mathbf{D}^b(\gr A)$. 

 (2) $\Leftrightarrow$ (3) It follows from Proposition \ref{finite-dim alg}.

(3) $\Rightarrow$ (4), (4) $\Rightarrow$ (1), and (4) $\Leftrightarrow$ (5) $\Leftrightarrow$ (6) are immediate.
\end{proof}

\begin{proposition} \label{relation between CM and Tor-p} Suppose that $A$ is a noetherian $\mathbb{N}$-graded algebra with $A_0$ semisimple and that $A$ has a balanced dualizing complex $R$. If $\xi=(\xi_0,\xi_1)\in \mathbb{R}_{\geqslant 0}\times \mathbb{R}$ with $\xi_1\leqslant \xi_0$, then, for any $0\neq X \in \mathbf{D}^b(\gr A)$ with $\pdim X< \infty$,
\begin{align} \label{stronger-inequality}
\Torreg_{\xi}(X) \leqslant \CMreg_{\xi}(X)+\Exreg_{-\xi}({}_AA),
\end{align}
and, for any $0\neq X \in \mathbf{D}^b(\gr A)$ with $\injdim X < \infty$,
\begin{align}\label{dual of stronger-inequality}
    \Exreg_{-\xi}(X) \leqslant -\ideg_{\xi}(X)+\Exreg_{-\xi}(A_A).
\end{align}   
In particular, $\CMreg_{\xi}(A)+\Exreg_{-\xi}({}_AA)\geqslant 0$ and $\CMreg_{\xi}(A)+\Exreg_{-\xi}(A_A)\geqslant 0$.
\end{proposition}

The proof follows the approach used in \cite[Proposition 5.4 and 5.9 (3)]{WY}. The inequality (\ref{dual of stronger-inequality}) follows directly from the duality between $\mathbf{D}_{fg}(\Gr A)$ and $\mathbf{D}_{fg}(\Gr A^o)$, together with the right version of the inequality (\ref{stronger-inequality}).

\begin{remark}\label{Torreg(R)=Torreg(A_0)}
    The inequality \eqref{stronger-inequality} is stronger than the inequality in Theorem \ref{relation between CM and Tor} (2) if $A_0$ is semisimple and $\xi=(\xi_0,\xi_1)\in \mathbb{R}_{\geqslant 0}\times \mathbb{R}$ with $\xi_1\leqslant \xi_0$. 
    If, furthermore, $A_0$ is semisimple and $\injdim {}_{A}A =\injdim A_A < \infty$, then $\pdim R_A < \infty$ and $\pdim {}_AR < \infty$. 
    The proofs of  Proposition \ref{relation between CM and Tor-p} and its right-handed version imply that
    $\Torreg_{\xi}(R_A)=\Torreg_{\xi}({}_{A}R)$ 
    and $\Exreg_{-\xi}({}_AA)=\Exreg_{-\xi}(A_A)$, where $\xi=(\xi_0,\xi_1)\in \mathbb{R}_{\geqslant 0}\times \mathbb{R}$ with $\xi_1\leqslant \xi_0$.
    
    If $A_0$ is semisimple, $\gldim A < \infty$ and $\xi=(\xi_0,\xi_1)\in \mathbb{R}_{\geqslant 0}\times \mathbb{R}$ with $\xi_1\leqslant \xi_0$, then, by Theorem \ref{relation between CM and Tor} (2) and Proposition \ref{relation between CM and Tor-p}, $\Torreg_{\xi}(R)=\Exreg_{-\xi}(A)=\Torreg_{\xi}(A_0)$.
    
\end{remark}

\subsection{Weighted CM-regularity homogeneous and Ex-regularity homogeneous} In this subsection, we introduce the weighted CM-regularity homogeneous and weighted Ex-regularity homogeneous properties for noetherian $\mathbb{N}$-graded algebras. We examine the conditions under which
the inequalities in Theorem \ref{relation between CM and Tor} become equalities, generalize \cite[Theorem 3.10]{KWZ2} to $\mathbb{N}$-graded algebras, and answer a question posed by Kirkman, Won, and Zhang in \cite[Remark 3.11 (2)]{KWZ2}.

\begin{definition}\label{CM-homogeneous-def} Let $A$ be an $\mathbb{N}$-graded algebra and $\xi=(\xi_0,\xi_1)\in \mathbb{R}\times \mathbb{R}$.
 \begin{itemize}
     \item [(1)] $A$ is called {\it  left (resp. right) $\xi$-CM-regularity homogeneous} if $\CMreg_{\xi}(A)=\CMreg_{\xi}(Ae_i)$  (resp. $\CMreg_{\xi}(A)=\CMreg_{\xi}(e_iA)$) for all $1 \leqslant i \leqslant n$.
     \item [(2)] $A$ is called {\it left (resp. right) $\xi$-Ex-regularity homogeneous}  if  $\Exreg_{\xi}({}_{A}A)=\Exreg_{\xi}(Ae_i)$ (resp. 
 $\Exreg_{\xi}(A_A)=\Exreg_{\xi}(e_iA)$) for all $1 \leqslant i \leqslant n$.
    \end{itemize}
    \end{definition}
    
\begin{remark}\label{ex-Tor homogeneous}
    If $A$ has a balanced dualizing complex $R$, then for all $1 \leqslant i \leqslant n$ and all $\xi=(\xi_0,\xi_1)\in \mathbb{R}\times \mathbb{R}$, $\Exreg_{\xi}({}_{A}A)=\Exreg_{\xi}(Ae_i)$ (resp. $\Exreg_{\xi}(A_A)=\Exreg_{\xi}(e_iA)$) if and only if $\Torreg_{-\xi}(R_A)=\Torreg_{-\xi}(e_iR)$ (resp. $\Torreg_{-\xi}({}_{A}R)=\Torreg_{-\xi}(Re_i)$).
\end{remark}

Every connected graded algebra $A$ is both $\xi$-CM-regularity homogeneous and $\xi$-Ex-regularity homogeneous for any $\xi=(\xi_0,\xi_1)\in \mathbb{R}\times \mathbb{R}$.  

Theorem \ref{relation between CM and Tor 2} below characterizes when the inequality in Theorem \ref{relation between CM and Tor} (2) is an equality. It extends \cite[Theorem 3.10]{KWZ2} and \cite[Theorem 5.24]{WY} and follows directly from Proposition \ref{sdeg-ideg equality 3} by taking $Y=R$, where $R$ denotes the balanced dualizing complex of $A$. Moreover, if $A$ is connected graded, then Theorem \ref{relation between CM and Tor 2} implies that, for any $\xi=(\xi_0,\xi_1)\in \mathbb{R}_{>0}\times \mathbb{R}$ with $\xi_1\leqslant \xi_0$, the equality  $\CMreg_{\xi}(X)=\Torreg_{\xi}(X)+\CMreg_{\xi}(A)$ holds for all $0\neq X\in \mathbf{D}^b(\gr A)$ with finite projective dimension. This answers a question posed by Kirkman-Won-Zhang in \cite[Remark 3.11 (2)]{KWZ2}. 

\begin{theorem}\label{relation between CM and Tor 2}
Let $A$ be a noetherian $\mathbb{N}$-graded algebra such that $A_0$ is semisimple and $A$ has a balanced dualizing complex $R$. Let $\xi=(\xi_0,\xi_1)\in \mathbb{R}_{>0}\times \mathbb{R}$ with $\xi_1\leqslant \xi_0$. Then the following statements are equivalent.
   \begin{itemize}
       \item [(1)] $A$ is left $\xi$-CM-regularity homogeneous.
       \item [(2)] $\CMreg_{\xi}(X)=\Torreg_{\xi}(X)+\CMreg_{\xi}(A)$ for all $0\neq X\in \mathbf{D}^b(\gr A)$ with finite projective dimension.
   \end{itemize}
\end{theorem}
\begin{proof}
By Definition \ref{CM-reg-def}, 
$$\CMreg_{\xi}(A) =-\ideg_{\xi}(\D(R\Gamma_A(A))) =-\ideg_{\xi}(R). $$
It follows from \cite[Theorem 5.1]{VdB} that, for all $1\leqslant i \leqslant n$,  
\begin{align*}
    \CMreg_{\xi}(Ae_i)
    =&-\ideg_{\xi}(\D(R\Gamma_A(Ae_i)))\\
    =&-\ideg_{\xi}(R\gHom_A(Ae_i,R))\\
    =&-\ideg_{\xi}(e_iR).
\end{align*}
Therefore, $\CMreg_{\xi}(A)=\CMreg_{\xi}(Ae_i)$ if and only if $-\ideg_{\xi}(R)=-\ideg_{\xi}(e_iR)$.

Similarly,
\begin{align*}
    \CMreg_{\xi}(X)
    =&-\ideg_{\xi}(\D(R\Gamma_A(X)))\\
    =&-\ideg_{\xi}(R\gHom_A(X,R)).
\end{align*} 
It follows from Proposition \ref{sdeg-ideg equality 3} that for any $X\in \mathbf{D}^b(\gr A)$ with finite projective dimension, 
$$\CMreg_{\xi}(X)=-\ideg_{\xi}(R\gHom_A(X,R))=\Extreg_{\xi}(X)-\ideg_{\xi}(R)$$
if and only if $\ideg_{\xi}(e_iR)=\ideg_{\xi}(R)$ for all $1\leqslant i\leqslant n$.
Hence,
$$\CMreg_{\xi}(X)=\Torreg_{\xi}(X)+\CMreg_{\xi}(A)$$ holds for all $0\neq X\in \mathbf{D}^b(\gr A)$ with finite projective dimension if and only if $A$ is left CM-regularity homogeneous.
\end{proof}

When $\xi=(1,0)$, Theorem \ref{relation between CM and Tor 2} yields the relation between the $a$-invariant and the $b$-invariant, as stated in the following corollary.

\begin{corollary}\label{relation of a-invariant}
    Suppose that $A$ is a noetherian $\mathbb{N}$-graded algebra such that $A_0$ is semisimple and $A$ has a balanced dualizing complex $R$. Then the following conditions are equivalent.
   \begin{itemize}
       \item [(1)] $A$ is left $(1,0)$-CM-regularity homogeneous.
       \item [(2)] $a^*(X)=b^*(X)+a^*(A)$ for all $0\neq X\in \mathbf{D}^b(\gr A)$ with finite projective dimension.
   \end{itemize}
\end{corollary}

\begin{proposition}\label{relation between CM and Tor 3}
    Suppose that $A$ is a noetherian $\mathbb{N}$-graded algebra such that $A_0$ is semisimple and $A$ has a balanced dualizing complex $R$. Let $\xi=(\xi_0,\xi_1)\in \mathbb{R}_{>0}\times \mathbb{R}$ with $\xi_1\leqslant \xi_0$. If $\injdim {}_{A}A=\injdim A_A <\infty$, then the following statements are equivalent.
    \begin{itemize}
        \item [(1)] $A$ is right $(-\xi)$-Ex-regularity homogeneous. 
        \item [(2)] $\Extreg_{\xi}(X)=\CMreg_{\xi}(X)+\Exreg_{-\xi}(A)$ for all $0\neq X\in \mathbf{D}^b(\gr A)$ with $\pdim X < \infty$.
       \item [(3)] $-\ideg_{\xi}(R\gHom_A(X,A))=\CMreg_{\xi}(X)+\Exreg_{-\xi}(A)$ for all $0\neq X\in \mathbf{D}^b(\gr A)$.
    \end{itemize}
    \end{proposition}
    \begin{proof}
  (1) $\Rightarrow$ (3) It follows from $R\gHom_{A^o}(R,e_iA)\cong R\gHom_A(\D(R\Gamma_{A^o}(e_iA)), A)\cong R\gHom_A(Re_i,A)$ that 
        \begin{align*}
            -\ideg_{\xi}(R\gHom_{A^o}(R,e_iA))=&-\ideg_{\xi}(R\gHom_A(Re_i,A)) \\
            =&\Extreg_{\xi}(Re_i) \quad (\textrm{by Theorem \ref{Extreg}})\\
            =&\Torreg_{\xi}(Re_i) \quad (\textrm{by Lemma \ref{Extreg-Torreg}})\\
            =&\Torreg_{\xi}(R) \quad (\textrm{by (1) and Remark } \ref{ex-Tor homogeneous})\\
            =&\Extreg_{\xi}(R) \quad (\textrm{by Lemma \ref{Extreg-Torreg}}).
        \end{align*}
Since $\pdim R_A < \infty$, 
it follows from the right version of Proposition \ref{sdeg-ideg equality 4} that for any $0\neq X \in \mathbf{D}^b(\gr A)$, 
\begin{align*}
    -\ideg_{\xi}(R\gHom_A(X,A))
    =&-\ideg_{\xi}(R\gHom_{A^o}(R,\D(R\Gamma_A(X))))\\
    =&\Extreg_{\xi}(R_A)-\ideg_{\xi}(\D(R\Gamma_A(X)))\\
    =&\Exreg_{-\xi}(A)+\CMreg_{\xi}(X).
\end{align*}

  (3) $\Rightarrow$ (2) If $\pdim X < \infty$, then by Theorem \ref{Extreg}, 
  $$\Torreg_{\xi}(X)=-\ideg_{\xi}(R\gHom_A(X,A))=\CMreg_{\xi}(X)+\Exreg_{-\xi}(A).$$

  (2) $\Rightarrow$ (1) Since $\D(R\Gamma_A(Re_i))\cong e_iA$, we have $\CMreg_{\xi}(Re_i)=0$ for any $1\leqslant i \leqslant n$. Taking $X=Re_i$ in (2) yields 
$$\Exreg_{-\xi}(e_iA)=\Extreg_{\xi}(Re_i)=\CMreg(Re_i) + \Exreg_{-\xi}(A)=\Exreg_{-\xi}(A).$$ 
 Thus, $A$ is right $(-\xi)$-Ex-regularity homogeneous. 
    \end{proof}

If $\gldim A < \infty$, then the following is a direct corollary of Proposition \ref{relation between CM and Tor 3}, in view of the equality $\Torreg_{\xi}(A_0)=\Exreg_{-\xi}(A)$ stated in Remark \ref{Torreg(R)=Torreg(A_0)}. 
   
\begin{corollary}\label{right ex-homo}
       Suppose that $A$ is a noetherian $\mathbb{N}$-graded algebra such that $A_0$ is semisimple and $A$ has a balanced dualizing complex $R$. Let $\xi=(\xi_0,\xi_1)\in \mathbb{R}_{>0}\times \mathbb{R}$ with $\xi_1\leqslant \xi_0$. If $\gldim A < \infty$, then the following statements are equivalent
       \begin{itemize}
        \item [(1)] $\Torreg_{\xi}(X)=\CMreg_{\xi}(X)+\Torreg_{\xi}(A_0)$ for all $0\neq X\in \mathbf{D}^b(\gr A)$.
       \item [(2)] $A$ is right $(-\xi)$-Ex-regularity homogeneous. 
    \end{itemize}
    \end{corollary}
 
\section{Relationship between weighted homological regularities and the AS-regular property}
In this section, following the ideas in\cite{KWZ2} and \cite{WY}, we study the relationship between weighted homological regularities and the $\mathbb{N}$-graded AS-regular property for noetherian $\mathbb{N}$-graded algebras, thereby generalizing results for connected graded algebras and the non-weighted case of $\mathbb{N}$-graded algebras.
 
A noetherian $\mathbb{N}$-graded algebra $A$ is said to satisfy the left {\it Auslander-Buchsbaum Formula} if, for any $0\neq X \in \mathbf{D}^b(\gr A)$ with $\pdim X < \infty$, $$\pdim X +\depth X =\depth A,$$ 
where $\depth{X}:=\mathop{\inf} \{ i\in \mathbb{Z}\mid \gExt^i_A (S,X)\neq 0 \}.$

Any noetherian connected graded algebra satisfying the $\chi$-condition satisfies the left Auslander-Buchsbaum formula \cite[Theorem 3.2]{Jo1}. For $\mathbb{N}$-graded algebras $A$ with a balanced dualizing complex, the following lemma characterizes when $A$ satisfies the Auslander-Buchsbaum formula.  

 \begin{theorem}\cite[Theorem 1.4]{LW}
     Let $A$ be a noetherian $\mathbb{N}$-graded algebra with a balanced dualizing complex. Then $A$ satisfies the left Auslander-Buchsbaum Formula if and only if $\gExt_{A^o}^d(M,A)\neq 0$ for all graded simple $A^o$-module $M$, where $d=\depth{A}$.
 \end{theorem}
 
Assuming the Auslander-Buchsbaum Formula holds, the following characterization of $\mathbb{N}$-graded AS-Gorenstein algebras extends the results of \cite[Theorem 3.6]{DW}.

\begin{proposition} \cite[Theorem 8.14]{LW} \label{ABF}
 Let $A$ be a noetherian $\mathbb{N}$-graded algebra with a balanced dualizing complex $R$. If $A$ satisfies the left and right Auslander-Buchsbaum Formula, then the following are equivalent.
  \begin{itemize}
      \item [(1)] $A$ is an $\mathbb{N}$-graded AS-Gorenstein algebra.
      \item [(2)] $\injdim({}_{A}A) < \infty$.
      \item [(3)] $\pdim({}_{A}R) < \infty$.
      \item [(4)] For any $X \in \mathbf{D}^b(\gr A)$, $\pdim(X) < \infty$ if and only if  $\injdim(X) < \infty$.
  \end{itemize}
\end{proposition}

We next use left weighted CM-regularities homogeneous property to establish sufficient conditions under which the Auslander-Buchsbaum formula holds for $\mathbb{N}$-graded algebras $A$ with a balanced dualizing complex.

\begin{lemma}\label{AB formula}
Suppose that $A$ is a noetherian $\mathbb{N}$-graded algebra with $A_0$ semisimple and a balanced dualizing complex $R$. Let $\xi=(\xi_0,\xi_1)\in \mathbb{R}_{>0}\times \mathbb{R}$. If $A$ is left $\xi$-CM-regularity homogeneous for $\xi_1 \ll 0$, then 
\begin{itemize}
    \item [(1)] $\depth {}_{A}A=\pdim X +\depth X$ for all $0 \neq X \in \mathbf{D}^b(\gr A)$ with finite projective dimension.
    \item [(2)] For all $0 \neq X \in \mathbf{D}^b(\gr A)$ with finite projective dimension,
    \begin{align*}
        &\sdeg (R^{\depth X}\Gamma_A(X)) \\
        =&\sdeg(\Tor^A_{\pdim X}(A_0,X))+\sdeg (R^{\depth A}\Gamma_A(A)).
    \end{align*}
\end{itemize}
\end{lemma}
\begin{proof}
     It follows from Theorem \ref{relation between CM and Tor 2} and $\xi_1\ll 0 \leqslant \xi_0$ that, for any $0 \neq X \in \mathbf{D}^b(\gr A)$ with finite projective dimension, $\CMreg_{\xi}(X)=\Torreg_{\xi}(X)+\CMreg_{\xi}(A)$. Since $\xi_1 \ll 0$,
    \[\CMreg_{\xi}(X)=\xi_1 \depth X + \xi_0 \sdeg (R^{\depth X}\Gamma_A(X)),\]
    \[\Torreg_{\xi}(X)=\xi_0 \sdeg(\Tor^A_{\pdim X}(A_0,X))-\xi_1 \pdim X.\]
In particular, 
    \[\CMreg_{\xi}(A)=\xi_1 \depth A +\xi_0 \sdeg (R^{\depth A}\Gamma_A(A)).\]
Hence, $\depth A= \depth X + \pdim X$ and 
  \[ \sdeg (R^{\depth X}\Gamma_A(X))=\sdeg(\Tor^A_{\pdim X}(A_0,X))+\sdeg (R^{\depth A}\Gamma_A(A)).\]
\end{proof}

\begin{lemma}\cite[Lemma 6.3]{WY}\label{graded equivalence}
 Let $A$ be an $\mathbb{N}$-graded algebra. Let $P=\mathop{\bigoplus}\limits_{i=1}^n Ae_{i}(p_i)$, where $p_1, p_2,\dots, p_n$ are integers, and let $B=\gEnd_A(P)$. Then 
 \begin{itemize}
     \item [(1)] $A$ satisfies the right Auslander-Buchsbaum
      formula if and only if $B$ does. 
     \item [(2)] $A$ has a balanced dualizing complex if and only if $B$ does.
     \item [(3)] \cite[Lemma 2.8]{IKU}\label{graded-Morita-equi} $A$ is a generalized AS-Gorenstein (regular) algebra if and only if $B$ is.
 \end{itemize}
\end{lemma}

To prove Theorem \ref{AGP}, the main result of this section, we adapt \cite[Theorem 4.7]{IKU} to ensure that the algebra $B$ appearing there is graded Morita equivalent to $A$ and that $B_0$ is semisimple. Note that \cite{IKU} considers right graded modules, whereas we consider left graded modules. The following theorem was proved in \cite[Theorem 6.4]{WY}.

 \begin{theorem}\label{AGP}\cite[Theorem 6.4]{WY}
     Let $A$ be a ring-indecomposable basic $\mathbb{N}$-graded AS-Gorenstein algebra of dimension $d$ with the average Gorenstein parameter $\ell^A_{av}\in\mathbb{Z}$. If $A_0$ is semisimple,
  then there exists a ring-indecomposable basic $\mathbb{N}$-graded AS-Gorenstein algebra $B$ of dimension $d$ whose Gorenstein parameters are all equal to $\ell^A_{av}$ such that 
  $B_0$ is semisimple and $B$ is graded Morita equivalent to $A$.
 \end{theorem}
   
For convenience, we make the following hypothesis.
 \begin{hypothesis}\label{Hy}
     $A$ is a noetherian (locally finite) $\mathbb{N}$-graded algebra with a balanced dualizing complex $R$, and $A_0$ is semisimple. 
 \end{hypothesis}
 
We are now ready to prove Theorems \ref{AS regular 11}, \ref{AS regular 22}, and \ref{AS regular 33} and their corollaries, thereby generalizing \cite[Theorems 3.2, 0.8]{KWZ1}. 
For any noetherian connected graded algebra $A$ with a balanced dualizing complex, \cite[Theorem 3.2, 0.8]{KWZ1} shows that the following conditions are equivalent: $A$ is AS-regular; $A$ is Cohen–Macaulay and $\ASreg(A) = 0$; $\ASreg(A) = 0$.

 \begin{theorem}\label{AS regular 11}
Suppose that $A$ satisfies Hypothesis \ref{Hy}, $A$ is basic and ring-indecomposable. Then the following conditions are equivalent.
\begin{itemize}
    \item [(1)] $A$ is $\mathbb{N}$-graded AS-regular of dimension $d$ with the average Gorenstein parameter $\ell_{av}^A \in \mathbb{Z}$.
    \item [(2)] For some integers $p_1, p_2,\dots, p_n$, $B: =\gEnd_A(\mathop{\bigoplus}\limits_{i=1}^n Ae_{i}(p_i))$ is an $\mathbb{N}$-graded CM-algebra of dimension $d$ with $B_0$ semisimple and $\ASreg(B)=0$.    
    \item [(3)] For some integers $p_1, p_2,\dots, p_n$, $B: =\gEnd_A(\mathop{\bigoplus}\limits_{i=1}^n Ae_{i}(p_i))$ is an $\mathbb{N}$-graded CM-algebra of dimension $d$ with $B_0$ semisimple and $\ASreg_{\xi}(B)=0$ for some $\xi=(\xi_0,\xi_1)\in \mathbb{R}_{>0}\times \mathbb{R}$.
\end{itemize}
\end{theorem}

\begin{proof}
(1) $\Rightarrow$ (2) See \cite[Theorem 6.6]{WY}. 
    
(2) $\Rightarrow$ (3) The conclusion follows from $\xi=(1,1)$ and (2).

(3) $\Rightarrow$ (1) 
Assume that $B_0=B_0e'_1\oplus B_0e'_2\oplus\cdots\oplus B_0e'_n$, where $\{e'_1,e'_2,\cdots,e'_n\}$ is a set of complete orthogonal primitive idempotents of $B_0$. 
It follows from $\ASreg_{\xi}(B)=0$ and Corollary \ref{ASreg=0} that $\CMreg_{\xi}(B)=\CMreg_{\xi}(Be'_i)$ for all $1\leqslant i\leqslant n$; hence, $B$ is left $\xi$-CM-regularity homogeneous. 

Since $B$ is a CM-algebra of dimension $d$, for any $1\leqslant i\leqslant n$,
\[\CMreg_{\xi}(B)=\xi_1 d+ \xi_0 \sdeg (R^d\Gamma_B(B)) \,\text{and} \]
\[\, \CMreg_{\xi}(Be'_i)=\xi_1 d+ \xi_0\sdeg (R^d\Gamma_A(Be'_i)).\]
Thus, for any $1\leqslant i\leqslant n$, $\sdeg (R^d\Gamma_B(B))=\sdeg (R^d\Gamma_B(Be'_i))$ and $d=\depth{B}=\depth{Be'_i}$.

It follows that, for any $\eta=(\eta_0,\eta_1) \in \mathbb{R}_{>0}\times \mathbb{R}$ and any $1\leqslant i\leqslant n$, $\CMreg_{\eta}(B)= \CMreg_{\eta}(Be'_i)$; that is, $B$ is left $\eta$-CM-regularity homogeneous. By Lemma \ref{AB formula}, $B$ satisfies the left Auslander-Buchsbaum formula. Similarly, the right-hand version shows that $B$ satisfies the right Auslander-Buchsbaum formula.

Suppose first that $\gldim B =d_1< \infty$. It follows from Lemma \ref{graded equivalence} that $B$ has a balanced dualizing complex. By Proposition \ref{ABF}, $B$ is an $\mathbb{N}$-graded AS-regular algebra of dimension $d_1$, say, with Gorenstein parameters $\{\ell_1,\ell_2,\cdots,\ell_n\}$. 
Since $B$ satisfies the left Auslander-Buchsbaum Formula, $$d_1=\injdim{B_B}=\pdim{{}_BR}=\depth{B}-\depth{R}=\depth{B}=d.$$
Hence 
\begin{align*}
    \CMreg_{\xi}(B)=&\xi_1 d+\xi_0 \sdeg (R^d\Gamma_B(B))\\
    =& \xi_1 d+\xi_0 \sdeg (\mathop{\bigoplus}\limits_{i=1}^n (\D(Be'_{\sigma (i)})(\ell_i))^{r_i}) \quad (\textrm {By \eqref{GAS-Goren} })\\
    =&\xi_1 d-\xi_0 \min\{\ell_i\}.
\end{align*}
Since $\ASreg_{\xi}(B)=0$,
$\Torreg_{\xi}(B_0) =- \CMreg_{\xi}(B) = - \xi_1 d + \xi_0 \min\{\ell_i\}.$

On the other hand, let $P^{\bullet}$ be the minimal graded projective resolution of ${}_{B}B_0$. It follows from the AS-regular property of $B$ that $P^{-d}=\mathop{\bigoplus}\limits_{i=1}^n Be_{\sigma(i)}(-\ell_{i})$.
Hence, $\sdeg(B_0\otimes_B P^{-d})=\max\{\ell_i\}$ and 
\begin{align*}
    \Torreg_{\xi}(B_0) =& \sup\{-\xi_1 i + \xi_0 \sdeg \Tor^B_i(B_0, B_0)\}\\
    \geqslant &-\xi_1 d + \xi_0 \sdeg(B_0\otimes_B P^{-d})=-\xi_1 d +\xi_0 \max\{\ell_i\}.
\end{align*} 
It follows that $- \xi_1 d + \xi_0 \min\{\ell_i\} \geqslant -\xi_1 d + \xi_0\max\{\ell_i\}$. Since $\xi_0>0$, 
$\min\{\ell_i\}=\max\{\ell_i\}=\ell_{av}^B$.
By Theorem \ref{AGP}, 
 $\ell_{av}^A=\ell_{av}^B \in \mathbb{Z}$. Since $A$ is graded Morita equivalent to $B$,
 $A$ is $\mathbb{N}$-graded AS-regular of dimension $d$ with Gorenstein parameter $\ell^A_{av}\in \mathbb{Z}$.

To complete the proof, it suffices to show that $\gldim B < \infty$. Suppose that $\gldim B = \infty$. Then $\pdim (B_0)=+\infty$. Let
 \begin{align*}
    \cdots\to P^{-j} \mathop{\to}^{\partial^{-j}} P^{-j+1}\to \cdots \to P^{-1} \to B \to B_0 \to 0
 \end{align*}
 be a minimal graded projective resolution of $B_0$ as a left $B$-module.  Then
 \begin{align*}
    \cdots\to P^{-j-2} \to P^{-j-1} \mathop{\to}^{\partial^{-j-1}} \ker \partial^{-j} \to 0
 \end{align*}
 is a minimal graded projective resolution, and $\pdim (\ker \partial^{-j})=+\infty$ for all $j \geqslant 0$.

 Since $B$ is a CM-algebra of dimension $d$, $\D(R\Gamma_B(B))\cong \D(R^d\Gamma_B(B))[d]:=\omega[d]$. 
 We claim that $\injdim_B \omega=d$ and $d=\cd(\Gamma_B)$. 
 
By \cite[Theorem 5.1]{VdB}, $\D(R\Gamma_B(X)) \cong \RHom_B(X,\omega[d])$ for all $X \in \mathbf{D}^b(\gr B)$.
Hence, for any $i> 0$ and $M \in \gr B$, $\gExt^i_B(M,\omega[d])\cong \D(R^{-i}\Gamma_B(M))=0$. Moreover, for any nonzero torsion module $N$, $\gExt^0_B(N,\omega[d])\cong \D(R^0\Gamma_B(N))\cong \D(N) \neq 0$. Consequently, $\injdim(\omega[d])=0$, and thus $\injdim \omega=d$. 
 
Since 
 $\D(R^q\Gamma_B(M))\cong \gExt^{-q}_B(M,\omega[d]) = \gExt^{-q+d}_B(X,\omega)=0$ for any $q> d$ and $M \in \Gr B$, it follows that $d=\cd(\Gamma_B)$.

Since $\RHom_B(B_0,\omega) =\RHom_B(B_0,\omega[d])[-d] \cong \D(R\Gamma_B(B_0))[-d]\cong B_0[-d]$,
 \begin{align*}
    \D(R^i\Gamma_B(\ker \partial^{-j} ))&\cong \gExt^{-i}_B(\ker \partial^{-j}, \omega[d])\\
    &\cong \gExt^{d-i}_B(\ker \partial^{-j},\omega)\\
    &\cong \gExt^{d-i+j+1}_B(B_0,\omega).
 \end{align*}
 Consequently, $\D(R^i\Gamma_B(\ker \partial^{-j} ))= 0$ for any $i< d$ and $j \geqslant d-1$.
Therefore, for all $j \geqslant d-1$, $\CMreg_{\xi}(\ker \partial^{-j})= \xi_1d +\xi_0 \sdeg(R^d\Gamma_B(\ker \partial^{-j}))$. 

If $j\geqslant d$, then the short exact sequence
    $0\to \ker \partial^{-j} \to P^{-j} \mathop{\to}^{\partial^{-j}} \ker \partial^{-j+1} \to 0$
 induces an exact sequence
\begin{align} \label{exact-seq-RGamma}
   0 \to R^d\Gamma_B(\ker \partial^{-j}) \to R^d\Gamma_B(P^{-j}) \to R^d\Gamma_B(\ker \partial^{-j+1}) \to 0.
\end{align}
Consequently,
\begin{equation}\label{CMreg-kernel}
\begin{aligned}
    \CMreg_{\xi}(\ker \partial^{-j})=&\xi_1 d+ \xi_0 \sdeg(R^d\Gamma_B(\ker \partial^{-j}))\\
    \leqslant &\xi_1 d +\xi_0 \sdeg( R^d\Gamma_B(P^{-j}))  \quad (\textrm{by } \eqref{exact-seq-RGamma})\\
    =&\CMreg_{\xi}(P^{-j})  \quad (B \textrm{ is Cohen-Macaulay})
\end{aligned}
\end{equation}
and 
\begin{equation}\label{sup-degree-B-tensor-P}
\begin{aligned}
    & \sdeg(B_0\otimes_B P^{-j-1})\\
    =&\xi_0^{-1}(-\xi_1 0 +\xi_0 \sdeg(B_0\otimes_B P^{-j-1}))
    \leqslant \xi_0^{-1} (\Torreg_{\xi}(\ker \partial^{-j}))\\
    \leqslant & \xi_0^{-1}(\CMreg_{\xi}(\ker \partial^{-j})+\Torreg_{\xi}(B_0))  \quad (\textrm{by Theorem} \ref{relation between CM and Tor} \,(2))\\
    \leqslant& \xi_0^{-1}(\CMreg_{\xi}(P^{-j})+\Torreg_{\xi}(B_0))  \quad (\textrm{by } \eqref{CMreg-kernel})\\
    =&\xi_0^{-1} \Torreg_{\xi}(P^{-j})  \quad (\textrm{by Corollary } \ref{ASreg=0})\\
    =& \sdeg(B_0\otimes_B P^{-j}).
\end{aligned}
\end{equation}
Let $c:=\Torreg_{\xi}(B_0)=-\CMreg_{\xi}(B)< \infty$, and let $c_1$ be a  positive integer such that $c_1>\xi_0^{-1}(c+(\xi_1-\xi_0)d)$. Since $\sdeg(B_0\otimes_B P^{-j})\geqslant j$ when $P^{-j}\neq 0$, it follows from \eqref{sup-degree-B-tensor-P} that
$$ d+c_1\leqslant \sdeg(B_0\otimes_B P^{-d-c_1})\leqslant \dots \leqslant \sdeg(B_0\otimes_B P^{-d}).$$
Hence, 
\begin{align*}
    \xi_0 \sdeg(B_0\otimes_B P^{-d})-\xi_1 d 
 &\geqslant \xi_0(d+c_1)-\xi_1 d \\
    &>  \xi_0(d+\xi_0^{-1}(c+(\xi_1-\xi_0)d))-\xi_1d \\
    &= \xi_0 d +c+(\xi_1-\xi_0)d-\xi_1 d=c,
\end{align*} 
which contradicts $c=\Torreg_{\xi}(B_0)$.
\end{proof}

Additional conditions equivalent to (1) in Theorem \ref{AS regular 11} follow from the characterization of $\ASreg_{\xi}(A)=0$ in Corollary \ref{ASreg=0}.

\begin{corollary}\label{coro ASreg=0}
  Suppose that $A$ satisfies Hypothesis \ref{Hy}, $A$ is basic and ring-indecomposable. Then the following are equivalent.
   \begin{itemize}
       \item [(1)] $A$ is $\mathbb{N}$-graded AS-regular of dimension $d$ with the average Gorenstein parameter $\ell_{av}^A \in \mathbb{Z}$.
       
   \item [(4)] $B=\gEnd_A(\mathop{\bigoplus}\limits_{i=1}^n Ae_{i}(p_i))$ is an $\mathbb{N}$-graded CM-algebra of dimension $d$ for some integers $ p_1, p_2,\cdots, p_n$, such that $B_0$ is semisimple and $\CMreg_{\xi}(X)-\Torreg_{\xi}(X)=c$ is constant for all $0 \neq X \in \mathbf{D}^b(\gr B)$ and for some $\xi=(\xi_0,\xi_1)\in\mathbb{R}_{>0}\times \mathbb{R}$.

   \item [(5)] $B=\gEnd_A(\mathop{\bigoplus}\limits_{i=1}^n Ae_{i}(p_i))$ is an $\mathbb{N}$-graded CM-algebra of dimension $d$ for some integers $ p_1, p_2,\cdots, p_n$, such that $B_0$ is semisimple and $\CMreg_{\xi}(M)-\Torreg_{\xi}(M)=c$ is constant for all $0 \neq M \in \gr B$ and for some $\xi=(\xi_0,\xi_1)\in\mathbb{R}_{>0}\times \mathbb{R}$.
   
    \item [(6)] $B=\gEnd_A(\mathop{\bigoplus}\limits_{i=1}^n Ae_{i}(p_i))$ is an $\mathbb{N}$-graded CM-algebra of dimension $d$ for some integers $ p_1, p_2,\cdots, p_n$, such that $B_0$ is semisimple and $\Exreg_{-\xi}(X)+\ideg_{\xi}(X)=-c$ holds for all $0 \neq X \in \mathbf{D}^b(\gr B)$ for some constant $c$ and some $\xi=(\xi_0,\xi_1)\in\mathbb{R}_{>0}\times \mathbb{R}$.
    
    \item [(7)] $B=\gEnd_A(\mathop{\bigoplus}\limits_{i=1}^n Ae_{i}(p_i))$ is an $\mathbb{N}$-graded CM-algebra of dimension $d$ for some integers $ p_1, p_2,\cdots, p_n$, such that $B_0$ is semisimple and $\Exreg_{-\xi}(M)+\ideg_{\xi}(M)=-c$ for all $0 \neq M \in \gr B$ and some $\xi=(\xi_0,\xi_1)\in\mathbb{R}_{>0}\times \mathbb{R}$.
   \end{itemize}
\end{corollary}

To remove the Cohen-Macaulay condition in part (3) of Theorem \ref{AS regular 11}, we introduce the following lemma.

\begin{lemma}\label{depth}
  Suppose that $A$ is an $\mathbb{N}$-graded algebra with finite local cohomology dimension. Let $d=\depth(A)$, $d_i=\depth(Ae_i)$ for all $1\leqslant i\leqslant n$, and $\xi=(\xi_0,\xi_1)\in \mathbb{R}_{>0}\times \mathbb{R}$. Then there exists $t \in \mathbb{R}$ such that for any $\xi_1\leqslant t\xi_0$, 
  \begin{align*}
      \CMreg_{\xi}(A)&=\xi_1 d +\xi_0\sdeg(R^d\Gamma_A(A)),\\
      \CMreg_{\xi}(Ae_i)&=\xi_1 d_i +\xi_0\sdeg(R^{d_i}\Gamma_A(Ae_i)).
  \end{align*}
\end{lemma}
\begin{proof}
    Let $m= \sup(R\Gamma_A(A))$, and $m_i=\sup (R\Gamma_A(Ae_i))$. It follows from 
    \begin{align*}
      \CMreg_{\xi}(A)&=\xi_1 d +\xi_0 \sdeg(R^d\Gamma_A(A)),\\
    \CMreg_{\xi}(Ae_i)&=\xi_1 d_i +\xi_0 \sdeg(R^{d_i}\Gamma_A(Ae_i)),
      \end{align*}
 that, for any $i\geqslant d$, 
  \begin{align*}
\xi_1 i +\xi_0 \sdeg(R^i\Gamma_A(A)) &\leqslant \xi_1 d +\xi_0\sdeg(R^d\Gamma_A(A)),\\
\xi_1 i +\xi_0 \sdeg(R^{i}\Gamma_A(Ae_i)) &\leqslant \xi_1 d_i +\xi_0 \sdeg(R^{d_i}\Gamma_A(Ae_i)).
\end{align*}
Let 
  \begin{align*}
      \gamma:&=\sdeg(R^d\Gamma_A(A))-\mathop{\max}\limits_{d\leqslant j\leqslant m}\{\sdeg(R^j\Gamma_A(A))\};\\
      \gamma_i:&= \sdeg(R^{d_i}\Gamma_A(Ae_i))-\mathop{\max}\limits_{d_i\leqslant j_i\leqslant m_i}\{ \sdeg(R^{j_i}\Gamma_A(Ae_i))\}.
  \end{align*}
  Then we may take $-\infty< t=\min \{ \gamma, \gamma_i \}\leqslant 0$, and the conclusions follow. 
\end{proof}

\begin{theorem}\label{AS regular 22}
 Suppose that $A$ satisfies Hypothesis \ref{Hy}, $A$ is basic and ring-indecomposable. Let $t$ be defined as in Lemma \ref{depth}. Then the following are equivalent.
     \begin{itemize}
         \item [(1)] $A$ is $\mathbb{N}$-graded AS-regular of dimension $d$ with average Gorenstein parameter $\ell_{av}^A \in \mathbb{Z}$.
         \item[(8)] $B=\gEnd_A(\mathop{\bigoplus}\limits_{i=1}^n Ae_{i}(p_i))$ is an $\mathbb{N}$-graded algebra for some integers $ p_1, p_2,\cdots, p_n$, such that $B_0$ is semisimple and $\ASreg_{\xi}(B)=0$ for some $\xi=(\xi_0,\xi_1)\in \mathbb{R}_{>0}\times \mathbb{R}$ with $\xi_1<t\xi_0$.
         \item[(9)] $B=\gEnd_A(\mathop{\bigoplus}\limits_{i=1}^n Ae_{i}(p_i))$ is an $\mathbb{N}$-graded algebra for some integers $ p_1, p_2,\cdots, p_n$, such that $B_0$ is semisimple and $\ASreg_{\xi}(B)=0$ for all $\xi=(\xi_0,\xi_1)\in \mathbb{R}_{>0}\times \mathbb{R}$ with $\xi_1\leqslant \xi_0$.
     \end{itemize}
\end{theorem}
 \begin{proof}
     (1) $\Rightarrow$ (9) By Theorem \ref{AGP} and Lemma \ref{graded-Morita-equi}, there exists a ring-indecomposable basic $\mathbb{N}$-graded AS-regular algebra $B=\gEnd_A(\mathop{\bigoplus}\limits_{i=1}^n Ae_{i}(p_i))$ of dimension $d$ with Gorenstein parameters $\{\ell_{av}^A,\ell_{av}^A,\cdots,\ell_{av}^A\}$ such that $B_0$ is semisimple. By Theorem \ref{def-GAS} and Definition \ref{GAS-Gorenstein}, $B$ has a balanced dualizing complex $R_B=B[d](-\ell_{av}^A)$. Consequently, for any $\xi=(\xi_0,\xi_1)\in \mathbb{R}_{>0}\times \mathbb{R}$ with $\xi_1\leqslant \xi_0$, $\CMreg_{\xi}(B)=-\ideg_{\xi}(R)=\xi_1 d-\xi_0 \ell_{av}^A$.

  By Proposition \ref{Extreg}, $\Extreg_{\xi}(X)=-\inf \{\xi_1i+\xi_0 j \mid \gExt^i_B(X,B)_j \neq 0 \}$ for any $X \in \mathbf{D}^b(\gr B)$.
    Thus 
    \begin{align*}
        \CMreg_{\xi}(X)=&-\ideg_{\xi}\big(D(R\Gamma(X))\big)=-\ideg_{\xi} \big(\RHom(X, R_B)\big)\\
        =&-\inf \{ \xi_1i+\xi_0 j \mid \gExt_B^i(X,R_B)_j \neq 0 \}\\
        =&-\inf \{\xi_1 i+\xi_0 j \mid \gExt_B^i(X, B[d](-\ell_{av}^A))_j \neq 0 \}\\
        =&\Extreg_{\xi}(X)+\xi_1 d-\xi_0 \ell_{av}^A=\Torreg_{\xi}(X)+\CMreg_{\xi}(B).
    \end{align*}
    It follows from Corollary \ref{ASreg=0} that $\ASreg_{\xi}(B)=0$ for all $\xi=(\xi_0,\xi_1)\in \mathbb{R}_{>0}\times \mathbb{R}$ with $\xi_1\leqslant \xi_0$.

     (9) $\Rightarrow$ (8) Obviously.

    (8) $\Rightarrow$ (1) By Theorem \ref{AS regular 11}, it suffices to show that $B$ is a CM-algebra of dimension $d$. 
        It follows from Corollary \ref{ASreg=0} that, for any $1\leqslant i \leqslant n$,  
    $ \CMreg_{\xi}(B)=\CMreg_{\xi}(Be'_i)$; that is, $B$ is left $\xi$-CM-regularity homogeneous. 
   By Lemma \ref{depth}, 
    \begin{align*}
      \CMreg_{\xi}(B)&=\xi_1 d +\xi_0 \sdeg(R^d\Gamma_B(B)),\\
      \CMreg_{\xi}(Be'_i)&=\xi_1 d_i +\xi_0 \sdeg(R^{d_i}\Gamma_B(Be'_i)),
  \end{align*}
   where $d=\depth{B}$ and $d_i=\depth{Be'_i}$. Suppose that there exists an $i$ such that $d_i > d$. Then $\xi_1 d > \xi_1 d_i$ by $\xi_1< 0$. Since $\xi_1 < t\xi_0$, 
  \begin{align*}
       &\xi_1 d_i + \xi_0 \sdeg(R^{d_i}\Gamma_B(B))-\xi_1 d-\xi_0 \sdeg(R^d\Gamma_B(B))\\
       =&\xi_1(d_i-d)-\xi_0 (\sdeg(R^d\Gamma_B(B))-\sdeg(R^{d_i}\Gamma_B(B)))\\
       <& \xi_0t(d_i-d)-\xi_0t\\
 =&\xi_0t(d_i-d-1)\leqslant 0.
  \end{align*}
  Thus, 
  \begin{align*}
      \xi_1 d+\xi_0 \sdeg(R^d\Gamma_B(B))&=\xi_1 d_i + \xi_0 \sdeg(R^{d_i}\Gamma_B(Be'_i)) \\
      &\leqslant \xi_1 d_i + \xi_0 \sdeg(R^{d_i}\Gamma_B(B))\\
      &< \xi_1 d+\xi_0 \sdeg(R^d\Gamma_B(B)),
  \end{align*}
 a contradiction. Hence, for any $1\leqslant i\leqslant n$, $d=d_i$ and  
\[\sdeg(R^d\Gamma_B(B))=\sdeg(R^{d_i}\Gamma_B(Be'_i))=\sdeg(R^{d}\Gamma_B(Be'_i)).\]
  
  It follows from Lemma \ref{depth} and $\xi_1 <\xi_0 t$ that, for any $\eta \leqslant \xi$ and any $1\leqslant i\leqslant n$, 
  \[ \CMreg_{\eta}(B)=\CMreg_{\eta}(Be'_i).\]
  Hence $B$ satisfies the left Auslander-Buchsbaum formula. A similar argument shows that $B$ satisfies the right Auslander-Buchsbaum formula.

 Suppose that $\gldim B < \infty$. By \cite[Lemma 3.6]{LW}, $\pdim {}_{B}R=\injdim (B_B) < \infty$. Let $P^{\bullet}$ be a minimal graded projective resolution of $R \in \mathbf{D}^b(\gr B)$. Then 
 \[ \pdim R = -\inf\{ r \in \mathbb{Z} \mid P^r\neq 0\} \geqslant -\inf R.\]
  Since $R\Gamma_B(R)\cong \D(B)$, $\depth R =0$ and $\depth B=\pdim R \geqslant -\inf R$. By the definition of $\depth B$,
  \[\depth B = \inf \{ r \in \mathbb{Z} \mid R^r\Gamma_B(B)\neq 0\} =- \sup R.\]
  Thus, $\sup R \leqslant \inf R$, and $B$ is a CM-algebra of dimension $d$.

  It suffices to prove $\gldim B < \infty$.
Suppose, for contradiction, that $\pdim B_0 = +\infty$. Since $\sdeg (\Tor_i^B(B_0, B_0)) \geqslant i$ for any $i$,   
\begin{equation*}
  \xi_0 \sdeg (\Tor_i^B(B_0, B_0))-\xi_1 i \geqslant \xi_0 i-\xi_1 i=(\xi_0-\xi_1)i.
\end{equation*}

 It follows from $\xi_1 < t\xi_0 < \xi_0$ that $\Torreg_{\xi}(B_0)=+\infty$. Hence $\CMreg_{\xi}(B)=-\infty$, contradicting $\CMreg_{\xi}(B)<\infty$. Thus $\gldim B < \infty$.
\end{proof}

The following theorem directly extends \cite[Theorem 0.6]{KWZ2} to $\mathbb{N}$-graded algebras. By \cite[Lemma 5.8]{WY}, for any $\mathbb{N}$-graded AS-Gorenstein algebra $A$, the condition $\Exreg(A)+\exreg(A)=0$ means that the Gorenstein parameters of $A$ are the same. This condition is also equivalent to $\Exreg_{\xi}(A)+\Exreg_{-\xi}(A)=0$ for some $\xi=(\xi_0,\xi_1)\in \mathbb{R}\times \mathbb{R}$ with $\xi_0\neq 0$.

\begin{lemma}\label{same gorenstein parameter}
    Let $A$ be an $\mathbb{N}$-graded AS-Gorenstein algebra of dimension $d$ with Gorenstein parameters $\{\ell_1,\ell_2,\cdots,\ell_n\}$. Then the following conditions are equivalent.
    \begin{itemize}
        \item [(1)] $\ell_1=\ell_2=\cdots=\ell_n$.
        \item [(2)]
        $\Exreg_{\xi}(A)+\Exreg_{-\xi}(A)=0$ for all $\xi=(\xi_0,\xi_1)\in \mathbb{R}\times \mathbb{R}$.
        \item [(3)] 
        $\Exreg_{\xi}(A)+\Exreg_{-\xi}(A)=0$ for some $\xi=(\xi_0,\xi_1)\in \mathbb{R}\times \mathbb{R}$ with $\xi_0\neq 0$.
     \end{itemize} 
\end{lemma}
\begin{proof} 
Since $A$ is an $\mathbb{N}$-graded AS-Gorenstein algebra of dimension $d$ with Gorenstein parameters $\{\ell_1,\ell_2,\cdots,\ell_n\}$, we have
$R\gHom_A(S,A)\cong \mathop{\bigoplus}\limits_{i=1}^n (e_{\sigma(i)}S(\ell_i))^{r_i} [-d]$.
Thus, for any $\xi_0\geqslant 0$,
$$\Exreg_{\xi}(A)=\xi_1d-\xi_0 \min_i\{\ell_i\}, \textrm{ and } \Exreg_{-\xi}(A)=-\xi_1d+\xi_0\max_i\{\ell_i\}.$$
Hence $\Exreg_{\xi}(A)+\Exreg_{-\xi}(A)=\xi_0(\max\{\ell_i\}-\min\{\ell_i\})$ for all $\xi_0\geqslant 0$.
Similarly, for any $\xi_0\leqslant 0$,
$$\Exreg_{\xi}(A)=\xi_1d-\xi_0 \max_i\{\ell_i\}, \textrm{ and } \Exreg_{-\xi}(A)=-\xi_1d+\xi_0\min_i\{\ell_i\}.$$
Consequently, $\Exreg_{\xi}(A)+\Exreg_{-\xi}(A)=-\xi_0(\max\{\ell_i\}-\min\{\ell_i\})$ for all $\xi_0\leqslant 0$.

(1) $\Rightarrow$ (2)
It follows from $\ell_1=\ell_2=\cdots=\ell_n$ that $\Exreg_{\xi}(A)+\Exreg_{-\xi}(A)=0$ for all $\xi=(\xi_0,\xi_1)\in \mathbb{R}\times \mathbb{R}$.

(2) $\Rightarrow$ (3) It is obvious.

(3) $\Rightarrow$ (1) 
It follows from $\xi_0\neq 0$ that $\Exreg_{\xi}(A)+\Exreg_{-\xi}(A)=0$ is equivalent to $\ell_1=\ell_2=\cdots=\ell_n$. 
\end{proof}

\begin{theorem}\label{AS regular 33} 
   Suppose that $A$ satisfies Hypothesis \ref{Hy}, and let $t$ be defined as in Lemma \ref{depth}. Then the following conditions are equivalent.
   \begin{itemize}
       \item [(1)] $A$ is $\mathbb{N}$-graded AS-regular of dimension $d$ and $\Exreg_{\xi}(A)+\Exreg_{-\xi}(A)=0$ for some $\xi=(\xi_0,\xi_1)\in \mathbb{R}\times \mathbb{R}$ with $\xi_0\neq 0$.
       
        \item [(2)] $\ASreg_{\xi}(A)=0$ for all $\xi=(\xi_0,\xi_1)\in \mathbb{R}_{>0}\times \mathbb{R}$ with $\xi_1\leqslant \xi_0$.
       \item [(3)] $\ASreg_{\xi}(A)=0$ for some $\xi=(\xi_0,\xi_1)\in \mathbb{R}_{>0}\times \mathbb{R}$ with $\xi_1< t\xi_0$. 
   \end{itemize}
    Furthermore, if $A$ is basic, then the following conditions are equivalent: $A$ is AS-regular over $A_0$; $\ASreg_{\xi}(A)=0$ for some $\xi_1 < \xi_0 t$; and $\ASreg_{\xi}(A)=0$ for all $\xi_1 \leqslant \xi_0$.
\end{theorem}

\begin{proof}
  (1) $\Rightarrow$ (2) 
By $\Exreg_{\xi}(A)=-\Exreg_{-\xi}(A)$ for some $\xi_0\neq 0$ and Lemma \ref{same gorenstein parameter}, $\ell := \max\{\ell_i\}=\min\{\ell_i\}$. By (\ref{GAS-Goren}), $ R\cong \mathop{\bigoplus}\limits_{i=1}^n (Ae_{\sigma(i)}(-\ell))^{r_i}[d]$, $\CMreg_{\xi}(A)=\xi_1d-\xi_0\ell$, and $\Torreg_{\xi}(R)=\xi_0\ell-\xi_1d$. Remark \ref{Torreg(R)=Torreg(A_0)} then implies that $\Torreg_{\xi}(A_0)=\Torreg_{\xi}(R)=\xi_0\ell-\xi_1d$. Hence $\ASreg_{\xi}(A)=0$ for all $\xi=(\xi_0,\xi_1)\in \mathbb{R}_{>0}\times \mathbb{R}$ with $\xi_1\leqslant \xi_0$.

(2) $\Rightarrow$ (3) Obviously.

(3) $\Rightarrow$ (1) Taking $B=A$ in the proof of $(8) \Rightarrow (1)$ in Theorem \ref{AS regular 22} shows that $A$ is $\mathbb{N}$-graded AS-regular of dimension $d$ with Gorenstein parameters $\{\ell,\ell,\cdots,\ell\}$. Hence, by Lemma \ref{same gorenstein parameter}, $\Exreg_{\xi}(A)+\Exreg_{-\xi}(A)=0$ for some $\xi=(\xi_0,\xi_1)\in \mathbb{R}\times \mathbb{R}$ with $\xi_0\neq 0$.

If $A$ is basic, it follows from Definition \ref{def-AS-Goren over A_0}  that $\ASreg(A)=0$ if and only if $A$ is AS-regular over $A_0.$
\end{proof}

If $A$ satisfies both the left and right Auslander–Buchsbaum formulas, then the following theorem is a direct generalization of \cite[Theorem 0.6]{KWZ2} and \cite[Theorem 6.10]{WY}.

\begin{theorem}\label{weighted AS regular 11}
    Suppose that $A$ satisfies Hypothesis \ref{Hy} and both the left and right Auslander-Buchsbaum formulas. Then the following conditions are equivalent. 
    \begin{itemize}
        \item [(1)]$A$ is $\mathbb{N}$-graded AS-regular of dimension $d$ and $\Exreg_{\xi}(A)+\Exreg_{-\xi}(A)=0$ for some $\xi=(\xi_0,\xi_1)\in \mathbb{R}\times \mathbb{R}$ with $\xi_0\neq 0$.
        
       \item [(2)] $\ASreg_{\xi}(A)=0$ for some $\xi=(\xi_0,\xi_1)\in \mathbb{R}_{>0}\times \mathbb{R}$ with $\xi_1\leqslant \xi_0$.
    \end{itemize}
       If $A$ is also basic, then $A$ is AS-regular over $A_0$ if and only if $\ASreg_{\xi}(A)=0$ for some $\xi=(\xi_0,\xi_1)\in \mathbb{R}\times \mathbb{R}$ with $\xi_1\leqslant \xi_0$.
\end{theorem}

\begin{proof}
  (1) $\Rightarrow$ (2) 
  It follows from (1) $\Rightarrow$ (2) in Theorem \ref{AS regular 33} that $\ASreg_{\xi}(A)=0$ for all $\xi=(\xi_0,\xi_1)\in \mathbb{R}_{>0}\times \mathbb{R}$ with $\xi_1\leqslant \xi_0$.

    (2) $\Rightarrow$ (1)
    If $\xi_1=\xi_0>0$, then, without loss of generality, we may assume $\xi_1=\xi_0=1$. The conclusion follows from \cite[Theorem 6.10]{WY}.

    Now, let $\xi_1<\xi_0$, and suppose that $\gldim A < \infty$. By Proposition \ref{ABF}, $A$ is $\mathbb{N}$-graded AS-regular of dimension $d$ with Gorenstein parameters $\{\ell_1, \ell_2,\cdots,\ell_n\}$. By (\ref{GAS-Goren}), $ R\cong \mathop{\bigoplus}\limits_{i=1}^n (Ae_{\sigma(i)}(-\ell_i))^{r_i}[d]$, $\CMreg_{\xi}(A)=\xi_1d-\xi_0\min\{\ell_i\}$ and $\Torreg_{\xi}(R)=\xi_0\max\{\ell_i\}-\xi_1d$. It follows from Remark \ref{Torreg(R)=Torreg(A_0)} that $\Torreg_{\xi}(A_0)=\Torreg_{\xi}(R)=\xi_0\max\{\ell_i\}-\xi_1d$. Hence, $0=\ASreg_{\xi}(A)=\CMreg_{\xi}(A)+\Torreg_{\xi}(A_0)=\xi_0(\max\{\ell_i\}-\min\{\ell_i\})$. Therefore, $\ell_1=\ell_2=\cdots=\ell_n$ and $\Exreg_{\xi}(A)+\Exreg_{-\xi}(A)=0$ by Lemma \ref{same gorenstein parameter}.
    
It remains to show that $\gldim A < \infty$.
Suppose, for contradiction, that $\pdim{A_0}=\infty$. Since $\sdeg (\Tor_i^A(A_0, A_0)) \geqslant i$ for all $i\geqslant 0$, 
 \[\xi_0\sdeg (\Tor_i^A(A_0, A_0))-\xi_1 i \geqslant \xi_0i-\xi_1 i=(\xi_0-\xi_1)i.\]
 Since $\xi_1<\xi_0$ and $\pdim{A_0}=\infty$, it follows that $\Torreg_{\xi}(A_0)\geqslant (\xi_0-\xi_1)\max\{i\mid i\in \mathbb{N}\}=\infty$. Hence $\CMreg_{\xi}(A)=-\Torreg_{\xi}(A_0)=\infty$, contradicting the finiteness of $\CMreg_{\xi}(A)$ established in Lemma \ref{CM-reg-finite}. Therefore, $\gldim A=\pdim{A_0} < \infty$.
\end{proof}

\subsection*{AI statement}
The authors used AI tools to edit the language, check the grammar, and improve the readability of the manuscript. The mathematical results and proofs were developed by the authors, who take full responsibility for the manuscript.


\thebibliography{plain}

\bibitem[AS]{AS} M. Artin and W. F. Schelter, Graded algebras of global dimension 3, Adv. Math. 66 (1987), 171--216.

\bibitem[AZ]{AZ} M. Artin and J.J. Zhang, Noncommutative projective schemes, Adv. Math. 109 (1994), 228–287.

\bibitem[BH]{BH} W. Bruns and J. Herzog, Cohen-Macaulay Rings, Cambridge Studies in Advanced Mathematics 39, Cambridge University Press, Cambridge, 1993.

\bibitem[CHT]{CHT} S. D. Cutkosky, J. Herzog and N. V. Trung, Asymptotic behaviour of Castelnuovo--Mumford regularity, Compositio Math. 118 (1999), 243--261.

\bibitem[DW]{DW}  Z.-C. Dong and Q.-S. Wu, Non-commutative Castelnuovo–Mumford regularity and AS regular algebras, J. Algebra 322 (2009), 122–136.

\bibitem[EG]{EG} D. Eisenbud and S. Goto, Linear free resolutions and minimal multiplicity, J. Algebra 88 (1984), 89–133.

\bibitem[GW]{GW} S. Goto and K. Watanabe, On graded rings, I, J. Math. Soc. Japan 30 Ž . 1978 ,179-213.

\bibitem[HT]{HT} L. T. Hoa and N. V. Trung, $a_i$-invariants and Castelnuovo--Mumford regularity of powers of ideals, J. Algebra 241 (2001), 154--178.

\bibitem[IKU]{IKU} O. Iyama, Y. Kimura and K. Ueyama, Cohen-Macaulay representations of Artin-Schelter Gorenstein algebras of dimension one, arxiv:2404.05925, 2024.

\bibitem[Jo1]{Jo1} P. J\o rgensen, Non-commutative graded homological identities, J. London Math. Soc. 57 (1998) 336–350.
 
\bibitem[Jo2]{Jo2} P. J\o rgensen, Non-commutative Castelnuovo-Mumford regularity, Math. Proc. Camb. Phil. Soc. 125 (1999), 203–221.

\bibitem[Jo3]{Jo3} P. Jørgensen, Noncommutative Castelnuovo--Mumford regularity for connected graded algebras, J. Algebra 211 (1999), 280--298.

\bibitem[Jo4]{Jo4} P. Jørgensen, Linear free resolutions over non-commutative algebras, Compositio Math. 140 (2004), 1053–1058.

\bibitem[JZ]{JZ} P. Jørgensen and J. J. Zhang, A gourmet's guide to Gorensteinness, Adv. Math. 151 (2000), 313--345.

\bibitem[KWZ1]{KWZ1} E. Kirkman, R. Won and J.J. Zhang, Homological regularities and concavities, Algebra and Number Theory, Vol. 20 (2026), No. 8, 1559--1595.

\bibitem[KWZ2]{KWZ2} E. Kirkman, R. Won and J.J. Zhang,
Weighted homological regularities, Trans. Amer. Math. Soc. 376 (2023), 7407--7445.

\bibitem[LW]{LW} H.-N Li and Q.-S Wu, Commonly graded algebras and their homological properties, arXiv: 2508.06068, 2025.

\bibitem[M]{M} D. Mumford, Lectures on curves on an algebraic surface, Princeton University Press, 1966.

\bibitem[MM]{MM} H. Minamoto and I. Mori, The structure of Artin-Schelter Gorenstein Algebras, Adv. Math. 226 (2011), 4061--4095.

\bibitem[MV1]{MV1} R. Martinez-Villa, Koszul algebras and the Gorenstein condition, Representations of algebras (S$\Tilde{\rm a}$o Paulo, 1999), 135–156, Lecture Notes in Pure and Appl. Math., 224, Dekker, New York, 2002.

\bibitem[MV2]{MV2} R. Martinez-Villa, Local Cohomology and non-commutative Gorenstein algebras, Int. J. Algebra 8 (2014), 329--348.

\bibitem[Ngu]{Ngu} H.D. Nguyen, Regularity bounds for complexes and their homology, Math. Proc. Camb. Phil. Soc.159 (2015), 355--377.

\bibitem[Ooi]{Ooi} A. Ooishi, Castelnuovo’s regularity of graded rings and modules, Hiroshima
Math. J. 12 (1982), 627--644. 

\bibitem[R\"{o}]{R\"{o}} T. R\"{o}mer, On the regularity over positively graded algebras, J. Algebra 319 (2008), 1--15.

\bibitem[RR]{RR} M. L. Reyes and D. Rogalski, Graded twisted Calabi-Yau algebras are generalized Artin-Schelter regular, Nagoya Math. J. 245 (2022), 100--153.

\bibitem[RRZ]{RRZ} M. L. Reyes, D. Rogalski and James J. Zhang, Skew Calabi–Yau algebras and homological identities, Adv. Math. 264 (2014), 308--354.

\bibitem[Sch]{Sch} P. Schenzel, Applications of Koszul homology to numbers of generators and syzygies, J. Pure Appl. Algebra 114 (1997), 287--303.

\bibitem[Tr1]{Tr1} N. V. Trung, Castelnuovo--Mumford regularity and degrees of generators of graded submodules, J. Algebra 200 (1998), 146--167.

\bibitem[Tr2]{Tr2} N. V. Trung. The largest non-vanishing degree of graded local cohomology modules, J. Algebra
215 (1999), 481--499.

\bibitem[Tr3]{Tr3} N.V. Trung, Castelnuovo-Mumford regularity and related invariants, Commutative algebra and combinatorics, Ramanujan Math. Soc. Lect. Notes Ser. 4 (2007), 157--180.

\bibitem[TW]{TW} N. V. Trung and H. J. Wang, On the Castelnuovo--Mumford regularity and the arithmetic degree, Math. Z. 249 (2005), 813--826.

\bibitem[VdB]{VdB} M. Van den Bergh, Existence theorems for dualizing complexes over noncommutative graded and filtered rings, J. Algebra 195 (1997), 662--679.

\bibitem[WY]{WY} Quanshui Wu, Bojuan Yi,  Numerical homological regularities over positively graded algebras, J. Algebra 686 (2026), 677--748.

\bibitem[Ye]{Ye} A. Yekutieli, Dualizing complexes over noncommutative graded algebras, J. Algebra 153 (1992), 41--84.

\bibitem[YZ]{YZ} A. Yekutieli, J.J. Zhang, Rings with Auslander dualizing complexes, J. Algebra 213 (1999), 1--51.

\bibitem[Z]{Z} J. J. Zhang, Dualizing complexes over noncommutative graded algebras, J. Algebra 218 (1999), 401--421.
\end{document}